\documentclass[11pt,reqno]{amsart} %put '[12pt,cram]' before {article}one can found
\usepackage{amsmath,amsthm,amsfonts,latexsym,euscript,amssymb,amscd}
\newtheorem{DF}{Definition}[section]
\newtheorem{LM}[DF]{Lemma}
\newtheorem{PROP}[DF]{Proposition}
\newtheorem{THM}[DF]{Theorem}
\newtheorem{COR}[DF]{Corollary}
\newtheorem{RMK}[DF]{Remark}
\newtheorem{RMKS}[DF]{Remarks}
\newtheorem{CONJ}[DF]{Conjecture}
\newtheorem{PROB}[DF]{Problem}
\newtheorem{QUES}[DF]{Question}
\newtheorem{EXM}[DF]{Example}

\newcommand{\bgdf}{\begin{DF}}
\newcommand{\nddf}{\end{DF}}
\newcommand{\bglm}{\begin{LM}}
\newcommand{\ndlm}{\end{LM}}
\newcommand{\bgprop}{\begin{PROP}}
\newcommand{\ndprop}{\end{PROP}}
\newcommand{\bgth}{\begin{THM}}
\newcommand{\ndth}{\end{THM}}
\newcommand{\bgthm}{\begin{THM}}
\newcommand{\ndthm}{\end{THM}}
\newcommand{\bgcor}{\begin{COR}}
\newcommand{\ndcor}{\end{COR}}
\newcommand{\bgrmk}{\begin{RMK}}
\newcommand{\ndrmk}{\end{RMK}}
\newcommand{\bgrmks}{\begin{RMKS}}
\newcommand{\ndrmks}{\end{RMKS}}
\newcommand{\bgconj}{\begin{CONJ}}
\newcommand{\ndconj}{\end{CONJ}}
\newcommand{\bgprob}{\begin{PROB}}
\newcommand{\ndprob}{\end{PROB}}
\newcommand{\bgques}{\begin{QUES}}
\newcommand{\ndques}{\end{QUES}}
\newcommand{\bgexm}{\begin{EXM}}
\newcommand{\ndexm}{\end{EXM}}
\newcommand{\bgequ}{\begin{equation}}
\newcommand{\ndequ}{\end{equation}}
\newcommand{\bgeq}{\begin{eqnarray}}
\newcommand{\ndeq}{\end{eqnarray}}
\newcommand{\bgeqq}{\begin{eqnarray*}}
\newcommand{\ndeqq}{\end{eqnarray*}}

\newcommand{\pf}{{\em Proof. }} %  \pf is already defined in old amslatex

\newcommand{\vv}{\vspace{2mm}\\}

\numberwithin{equation}{section}

\newcommand{\dfref}[1]{Definition~\ref{#1}}
\newcommand{\lemref}[1]{Lemma~\ref{#1}}
\newcommand{\lmref}[1]{Lemma~\ref{#1}}
\newcommand{\propref}[1]{Proposition~\ref{#1}}

\newcommand{\thmref}[1]{Theorem~\ref{#1}}

\newcommand{\quesref}[1]{Question~\ref{#1}}
\newcommand{\secref}[1]{\S\ref{#1}}
\title[Normal Quantum Subgroups, Properties $F$ and $FD$]
{\bf Equivalent Notions of Normal Quantum Subgroups,
Compact Quantum Groups with Properties $F$ and $FD$,  and Other Applications}
\author[Shuzhou Wang]
{\bf Shuzhou Wang}

\address{Department of Mathematics, University of Georgia,
Athens, GA 30602
\newline \indent
Fax: 706-542-2573; Tel: 706-542-0884
}
\email{szwang@math.uga.edu}

\begin{document}

\begin{abstract}
The notion of normal quantum subgroup introduced
in algebraic context by Parshall and Wang
when applied
to compact quantum groups is shown to be equivalent to
the notion of normal quantum subgroup introduced by the author.
As applications, a quantum analog of the
third fundamental isomorphism theorem for groups is obtained,
which is used along with the equivalence theorem to obtain results on
structure of quantum groups with property $F$ and quantum groups with property $FD$.
Other results on normal quantum subgroups for tensor products,
free products and crossed products are also proved.
\end{abstract}

\maketitle

\section{Introduction}

The notion of normal quantum subgroup is an important and subtle
concept in the theory of quantum groups.
In purely algebraic context of Hopf algebras,
B. Parshall and J. Wang \cite{ParshallWang91a} defined a
notion of normal quantum subgroup using left and right adjoint coactions of
the Hopf algebra on itself, which was further studied
by other authors such as Schneider \cite{Schneid93a}, Takeuchi  \cite{Takeuchi94a},
and Andruskiewitsch and Devoto \cite{Andrus95a}.
Parshall and Wang noted that left normal quantum groups may not be
right normal in general, and given a normal quantum subgroup in their sense,
it is not known whether there exists an associated exact sequence,
and if an exact sequence exists, it may not be unique.
These difficulties are peculiar phenomena of general Hopf algebras in purely algebraic
context distinguishing Hopf algebras from groups.
Other complications related to the notion of normal quantum groups
in purely algebraic context are included in \cite{Schneid93a}.
In $C^*$-algebraic context, the author introduced \cite{free}
a notion of normal quantum subgroup of compact quantum groups
 using analytical properties of representation theory of compact quantum groups.
 It was not known whether these two notions of normality are equivalent when
they are applied to the canonical dense Hopf $*$-algebras of quantum
representative functions of compact quantum groups. In  \cite{simple},
the author's notion of normal quantum groups was used
in an essential way to define the notion of simple compact quantum groups.
It was also announced in \cite{simple} without proof that the
above two notions of normality are equivalent for compact quantum groups
(see remark (b) after Lemma 4.4 in  \cite{simple}).
As consequences, left normal and right normal defined in algebraic context
by Parshall and Wang are also equivalent for compact quantum groups,
 and their normal quantum subgroups always
 give rise to a unique exact sequence. That is,  the
  complications mentioned above in purely algebraic setting
  do not present themselves in the world of compact quantum groups.
  Such properties might be useful for formulating an
  appropriate notion of quantum groups in algebraic setting, which
 is still an open problem.

  Other facts announced in \cite{simple} without proofs include general results on
structure of compact quantum groups with property $F$ (resp. property $FD$),
where, roughly speaking, a compact quantum group $G$ is said to have { property $F$}
if its quantum function algebra $A_G$ has the same property with respect to
quotients by normal quantum subgroups as the function algebra of a compact group,
and it is said to have { property $FD$}
if its quantum function algebra has the same property
with respect to quotients by normal quantum subgroups
as the quantum function algebra of the dual of a discrete group.
See \dfref{property-F-etc} below for precise definitions
of these concepts and notation used above.
Compact quantum groups with property $F$ include all quantum groups
obtained from compact Lie groups by deformation method, such as
compact real form of the Drinfeld-Jimbo quantum groups
and Rieffel's deformation, as well as most of the universal quantum groups
constructed by the author except the universal unitary quantum groups $A_u(Q)$
(also called the free unitary quantum groups), cf. \cite{simple}.

The purposes of this paper are
to give complete proof of the equivalence of the two notions of normality mentioned
above and give the following applications of this Equivalence Theorem on the
structure of compact quantum groups.

(1) We establish a complete quantum analog of the Third Fundamental Isomorphism Theorem.
This is the only one among the three fundamental isomorphism theorems that
has a complete quantum analog without added conditions or restrictions.
On the contrary, a surjection of compact quantum groups
(i.e. inclusion of Woronowicz $C^*$-algebras) does not always give rise to
a quantum analog of the First Fundamental Isomorphism
Theorem, except in the special case where an exact sequence can
be constructed, cf. \cite{ParshallWang91a,Schneid93a,Takeuchi94a,Andrus95a} for this and
other subtleties. Taking the example of the group $C^*$-algebra $A_G :=C^*(F_2)$
of the free group $F_2$ on two generators,
a Woronowicz $C^*$-subalgebra of $A_G$ does not give rise to an
exact sequence unless it is the group $C^*$-algebra of a normal subgroup
of $F_2$. In addition, it is not clear at the moment
how a quantum analog of the second fundamental isomorphism theorem
can be formulated.

(2) Using the Equivalence Theorem and the
quantum analog of the Third Fundamental Isomorphism Theorem,
we show that quotient quantum groups of a compact quantum group
with property $F$ also have property $F$,
and quantum subgroups of a compact quantum group
with property $FD$ also have property $FD$. We show that
quotient quantum groups of a compact quantum group
with property $FD$ also have property $FD$ provided
$G$ has the pullback property.
The pullback property is the quantum group version of the group situation
in which every subgroup of $G/N$ is of the form $H/N$ for some
subgroup $H$ of $G$ containing $N$.
We give an example to show not all compact quantum groups have
the pullback  property.

(3) We prove results on normal quantum subgroups for tensor products,
free products and crossed products.
Note that the free product construction has no place in the classical world
of compact groups. It is a total quantum phenomenon.

An outline of the paper is as follows.
In Section 2, we recall the algebraic notion of normal quantum subgroups in
\cite{ParshallWang91a} and the analytical notion of normal quantum subgroups
in \cite{free} respectively.
In Section 3, the equivalence of these two notions of normality is
proved.
In Section 4, as applications of the Equivalence Theorem,
we prove the quantum analog of the Third Fundamental Isomorphism Theorem,
and results on structure of compact quantum groups
with property $F$ and property $FD$. In Section 5, as further applications,
properties of normal quantum subgroups for free
products, tensor products and crossed products are given.

We note that most results in this paper, such as
\thmref{thm-equivalence} and those in Sections 4 and 5,
are also valid for cosemisimple Hopf algebras when they are
appropriately re-formulated.
For instance, one simply replaces the statement in (1)  of \thmref{thm-equivalence}
with the equality in (3)$^\prime$ of \propref{normal-subgroup} for
such a reformulation. The existence of the Haar integral/measure
shared by both compact quantum groups and cosemisimple Hopf algebras
is a key element in the proofs of these results.

Besides the general abstract theory on compact quantum groups
developed by Woronowicz and general constructions of particular classes of
compact quantum groups, there seem to be few general results on
the {\em structure} of {\em infinite} compact quantum groups in the
literature with the possible exception of \cite{CDPS12pr},
a situation contrary to finite quantum groups for
which there is much literature on their structure and classification.
The results in sections 4 and 5 are a modesty attempt at
developing theory on structure of infinite compact quantum groups.
It is expected that such results will be useful in the program  \cite{simple}
of classification of simple compact quantum groups
and further study of the structure of compact quantum groups.
\vv
{\bf Convention:} We use the notation and terminology in \cite{free}.
For a compact quantum group $G$, $A_G$ denotes the
underlying Woronowicz $C^*$-algebra and ${\mathcal A}_G$ the associated
canonical dense Hopf $*$-algebra of quantum representative functions on $G$.
Sometimes it is convenient to abuse the notation by calling $A_G$ a
compact quantum group, referring to $G$. As was pointed out on
p.533 of \cite{krein},
morphisms between quantum groups are
meaningful only for {\em full} Woronowicz $C^*$-algebras $A_G$
(i.e. restriction of the norm $||\cdot ||$ of
$A_G$ to the $*$-algebra ${\mathcal A}_G $
is the maximum of all possible $C^*$-norm on ${\mathcal A}_G $), although one can
define morphisms between arbitrary Woronowicz $C^*$-algebras (cf. 2.3 in \cite{free}).
Unless otherwise explicitly stated, we assume that {\em all} Woronowicz $C^*$-algebras
considered in this paper to be full. We also use standard notation in Hopf algebras, including
Sweedler's summation convention \cite{Sweedler,Montgomery93a}, and $\Delta$,
$\varepsilon$, $S$ for coproduct, counit and antipode, respectively.
Relevant basic information on compact quantum groups and Hopf algebras can also be
found in \cite{KlimykSchmudgen}.

\section{Two Notions of Normal Quantum Subgroups and Their Equivalence}
\label{notions-of-normal}

We recall the two notions of normal quantum subgroups defined by the author
in \cite{free} analytically and by Parshall and Wang in \cite{ParshallWang91a}
algebraically.

\bgdf
{\rm (cf. 2.3 and 2.13 in \cite{free})}
\label{qsg}
A {\bf quantum subgroup} of a compact quantum group
$G$ in the sense of \cite{free}
is a pair $(N, \pi)$, where $A_N$ is a Woronowicz $C^*$-algebra
and $\pi: A_G \longrightarrow A_N$ is a surjection of $C^*$-algebras
that satisfies
\bgeq
(\pi \otimes \pi) \Delta_G = \Delta_N \pi,
\ndeq
where $\Delta_G$ and $\Delta_N$ are the coproducts of $A_G$ and
$A_N$ respectively.
\nddf

It can be shown (cf. 2.9 and 2.11 in \cite{free} or 1.3.9 and 1.3.9 in \cite{Wang})
that $(N, \pi)$ is a quantum subgroup of $G$
if and only if the kernel $\ker(\pi)$ (denoted by ${I}$)
of $\pi$ is a {\bf Woronowicz $C^*$-ideal} of ${A}_G$ in the sense that it is a
$C^*$-ideal of  ${A}_G$ that satisfies
\bgeq
\label{coproduct0}
\Delta_G({I}) \subset \ker (\pi \otimes \pi).
\ndeq

When there is possible confusion, such as when referring to kernel
of the morphism in \lmref{reconstruct-N-from-G/N} and when comparing
analytically defined quantum subgroup with its algebraically defined
counterpart, we use
$\hat{\pi}: {\mathcal A}_G \longrightarrow {\mathcal A}_N$
to denote the restriction of $\pi$ that maps the dense algebra
${\mathcal A}_G $  of quantum representative functions on $G$ onto
that $ {\mathcal A}_N $ on $N$.
The quantum group $(N, \pi)$ should be more precisely called a closed
quantum subgroup, but we will omit the word {\em closed }
since we do not consider non-closed quantum subgroups.

For convenience of readers who are familiar with the language
of comodules of Hopf algebras but less so with the notion of
a finite dimensional representation of a compact quantum group $G$,
we recall the definition of the latter.
For definition of infinite dimensional (unitary) representations,
see 2.4 in \cite{free} or Section 3 in \cite{Wor98a}.

\bgdf
{\rm (cf. 2.1 in \cite{Wor87b})}
\label{rep}
A representation of dimensional $d$ of a compact quantum group $G$
 is an invertible element $v$ of the algebra $M_d(A_G)$ of $d \times d$
 matrices with entries in $A_G$ such that
\bgeq
 \Delta_G (v_{ij}) = \sum_{k=1}^d v_{ik} \otimes v_{kj}.
\ndeq
\nddf

As shown in Proposition 13 on p30 in
\cite{KlimykSchmudgen} {\em and} Proposition 3.2 in \cite{Wor87b},
finite dimensional representations of $G$ and
comodules of the Hopf algebra ${\mathcal A}_G$
are in natural one to one correspondence.
Note that representation in the sense above is called
non-degenerate representation in \cite{Wor87b}.

The algebra $M_d(A_G)$  is also written as $M_d({\mathbb C}) \otimes A$,
where $M_d({\mathbb C})$ is the algebra of $d \times d$
 matrices with entries in complex numbers  ${\mathbb C}$.
For this reason, for a linear map $\tau$ from $A_G$ to another vector space,
the matrix $(\tau(v_{ij}))_{i,j=1}^d$ is often written in three different ways
interchangeably with slight abuse of notation:
\bgeqq
(\tau(v_{ij}))_{i,j=1}^d = (id \otimes \tau) (v)  = \tau(v).
\ndeqq

\bgdf
{\rm (cf. 2.13 in \cite{free})}
\label{normal}
A quantum subgroup $(N, \pi)$ of $G$ is called {\bf normal} if for every
irreducible representation $u^\lambda$ of $G$,
the multiplicity of the trivial
representation of $N$ in $ \pi (u^\lambda)$ is equal to
either zero or the dimension of $u^\lambda$.
\nddf

Let $h_N$ be the Haar measure (also called Haar state or Haar integral) on $N$. Then it is clear that
$N$ is normal if and only if for every irreducible
representation $u^\lambda$ of $G$,
\bgeqq
\text{either} \; \;  h_N \pi (u^\lambda) = I_{d_\lambda}, \; \; \text{or} \; \;
h_N \pi (u^\lambda) = 0,
\ndeqq
where $d_\lambda$ is the dimension of $u^\lambda$, and
$I_{d_\lambda}$ is the $d_\lambda \times d_\lambda$ identity matrix.

We recall the definition of normal quantum subgroup in Parshall and
Wang \cite{ParshallWang91a} adapted to Hopf
*-algebras of the form ${\mathcal A}_G$ where $G$ is a compact quantum
group, though their definition applies to more general Hopf algebras.

\bgdf
{\rm (cf. 1.4 in \cite{ParshallWang91a})}
\label{aqsg}
An {\bf algebraic quantum subgroup} of a compact group $G$
is a pair $(N, \eta)$ where $N$ is a compact quantum $N$
and $\eta: {\mathcal A}_G \longrightarrow {\mathcal A}_N$
is a surjection of $*$-algebras that satisfies
\bgeq
\label{coproduct1}
(\eta \otimes \eta) \Delta_G = \Delta_N \eta ,
\; \; \\
\label{counit-antipode1}
\varepsilon_G = \varepsilon_N \eta ,
\; \;
\eta S_G = S_N \eta,
\ndeq
where $\eta_G$ and $\eta_N$ (resp. $S_G$ and $S_N$) are the counits (resp.  antipodes) of
${\mathcal A}_G$ and ${\mathcal A}_N$ respectively.
\nddf

It is clear that $(N, \eta)$ is an {\bf algebraic quantum subgroup} of $G$
if and only if the kernel $\ker(\eta)$ (denoted by ${\mathcal I}$)
of $\eta$ is a Hopf $*$-ideal of ${\mathcal A}_G$ in the sense that it is a
$*$-ideal of  ${\mathcal A}_G$ that satisfies
(cf.  1.4 in \cite{ParshallWang91a})
\bgeq
\label{coproduct2}
\Delta_G({\mathcal I}) \subset {\mathcal A}_G \otimes {\mathcal I}
+ {\mathcal I} \otimes {\mathcal A}_G,
\; \; \\
\label{counit-antipode2}
\varepsilon_G({\mathcal I}) = 0,
\; \;
S_G({\mathcal I}) \subset {\mathcal I} .
\ndeq

In \cite{ParshallWang91a}, the morphism $\eta$ is not
required to preserve the $*$-algebra structure and the ideal ${\mathcal I}$ is
not required to be a $*$-ideal. Since we restrict attention to
compact quantum groups,  we need to require both.

Using Woronowicz's Peter-Weyl theory for compact quantum groups,
one can easily show (cf. 2.10 of \cite{free} and details in 1.2.16 of \cite{Wang})
that if $(N, \pi)$ is a quantum subgroup of $G$ in the
sense of \dfref{qsg}, then $(N, \hat{\pi})$ is an algebraic quantum subgroup of $G$.
In particular, the counits and antipodes of the associated
dense Hopf subalgebras are {\em automatically} preserved, i.e. conditions in
\eqref{counit-antipode1} (resp. \eqref{counit-antipode2}) {\em automatically}
follow from \eqref{coproduct1} (resp. \eqref{coproduct2}), which are postulated
in 1.4 of Parshall and Wang \cite{ParshallWang91a}.
However it must be cautioned that if $(N, \eta)$ is an
algebraic quantum subgroup of $G$ in the original sense of \cite{ParshallWang91a}
without $\eta$ preserving the $*$-structures of ${\mathcal A}_G$ and ${\mathcal A}_N$
and both $G$ and $N$ are compact, there is {\em no} morphism of compact quantum
groups $\pi: A_G \longrightarrow A_N$ with $\hat{\pi} = \eta$.

The precise correspondence between analytical quantum subgroups in \dfref{qsg}
and algebraic quantum subgroups of $G$ in \dfref{aqsg} is given by the
following theorem (see 4.3.(2) in \cite{simple}), which essentially says that
the two notions are equivalent. It is the first
step that reduces the $C^*$-setting to algebraic setting for the proof of the
equivalence theorem on normality:

\bgthm
\label{Hopf $*$-ideals vs. Woronowicz $C^*$-ideals}

The map $f({\mathcal I})=\overline{\mathcal I}$ is a bijection
from the set of Hopf $*$-ideals $\{ {\mathcal I} \}$
of ${\mathcal A}_G$ onto the set of Woronowicz
$C^*$-ideals $\{ I \}$ of $A_G$ with full quotient Woronowicz
$C^*$-algebra $A_G / I$.
The inverse $g$ of $f$ is given by
%$f^{-1}(I)
$g(I) = I \cap {\mathcal A}_G$.

\ndthm

\noindent
{\em Remarks:} A detailed proof of the above theorem is given in \cite{simple}. We note that its
proof is a nice interplay between the algebraic and analytical properties
of compact quantum groups.  We also note that
many concrete constructions in the analytical $C^*$-algebraic context
 also have purely algebraic formulation, such as the quantum permutation groups in
 \cite{qsym} and their algebraic counter part in Bichon \cite{Bi03a, Bi08a}.
 The theory of compact quantum groups is a rich ground
 where algebraic aspects and analytical
 aspects pleasantly interplay with each other.

\medskip

Thanks to  \thmref{Hopf $*$-ideals vs. Woronowicz $C^*$-ideals}, we can
now focus on the algebraic object ${\mathcal A}_G$.
Let $a \in {\mathcal A}_G$. The left and right adjoint coactions are defined
respectively by
$$
ad_l(a) : = \sum a_{(2)} \otimes  a_{(1)} S(a_{(3)}),
\; \; \;  ad_r(a) : = \sum a_{(2)} \otimes S(a_{(1)}) a_{(3)} ,
$$
where $S$ is the antipode of the Hopf algebra ${\mathcal A}_G$ and
 Sweedler's notation \cite{Sweedler} is used:
$$(\Delta \otimes id) \Delta(a)= (id \otimes \Delta) \Delta(a)
=\sum a_{(1)} \otimes a_{(2)} \otimes a_{(3)}. $$

\bgdf
{\rm (cf. Definition 1.5 in \cite{ParshallWang91a})}
\label{a-normal}
An algebraic quantum subgroup $(N, \eta)$ of $G$ is called {\bf a-normal} if
 $\ker(\eta)$
is a normal Hopf ideal of ${\mathcal A}_G$ in the sense
that the following two conditions are satisfied,
$$
%i.e, \; \;
ad_l(a) \in \ker(\eta) \otimes {\mathcal A}_G,
\; \; \; {\text and}
\; \; \;
ad_r(a) \in \ker(\eta) \otimes {\mathcal A}_G
$$
for all $a \in \ker(\eta)$.
\nddf

To compare with our \dfref{normal},
for the time being we use the term {\bf a-normal} instead of
{normal} for the situation considered by Parshall and Wang
\cite{ParshallWang91a}. Following their paper
we call $(N, \eta)$ {\bf left a-normal} (resp. {\bf right a-normal}) if the first
(resp. second) condition in \dfref{a-normal} above is satisfied.
In Schneider \cite{Schneid93a}, a morphism such as $\eta$ used in
\dfref{a-normal} above is also called a {conormal} morphism. In 1.1.7 of
Andruskiewitsch and Devoto \cite{Andrus95a},  ${\mathcal A}_N$ is said be a
{right quotient ${\mathcal A}_G$ comodule} if the second condition above is satisfied,
because the comodule structure $ad_r$ can then be induced to the quotient
${\mathcal A}_G/\ker(\eta)$, which is ${\mathcal A}_N$.

Our {\em first goal} in this paper is to prove
the following Equivalence Theorem.

\bgthm
\label{thm-equivalence}
Let $(N, \pi)$ be a quantum subgroup of a compact quantum group $G$.
Then the following conditions are equivalent.

{\rm (1)}  $(N, \pi)$ is normal.

{\rm (2)}  $(N, \hat{\pi})$ is a-normal.

{\rm (3)}  $(N, \hat{\pi})$ is left a-normal.

{\rm (4)}  $(N, \hat{\pi})$ is right a-normal.
\ndthm

The proof is given in the next section.

\section{Proof of \thmref{thm-equivalence}}
\label{proof-equivalence}

For convenience of the reader, we recall the notations to be used below.
Define
$$
A_{G/N} = \{ a \in A_G | (id \otimes \pi) \Delta (a) = a \otimes 1_N \},
\; \; \;
$$
$$
A_{N \backslash G} = \{ a \in A_G | (\pi \otimes id) \Delta (a)
= 1_N \otimes a \}, \; \; \;
$$
where $\Delta$ is the coproduct on $A_G$, $id$ is
the identity map on $A_G$, and $1_N$ is the unit of the algebra
$A_N$, which will simply be denoted by $1$ when the context is clear.
Similarly, we define
$$
{\mathcal A}_{G/N} = {\mathcal A}_G \cap
A_{G/N}, \; \; \;
\text{and}
\; \; \;
{\mathcal A}_{N \backslash G} =
{\mathcal A}_G \cap
A_{N \backslash G} .
$$
Note that $G/N$ and $N \backslash G$ should be denoted more precisely
by $G/(N, \pi)$ and $(N, \pi) \backslash G$ respectively if there is
a possible confusion. Let $h_N$ be the Haar measure on $N$.
Let
$$
E_{G/N} = ( id \otimes h_N \pi ) \Delta
, \; \; \;
\; \; \;
E_{N \backslash G} = (h_N \pi \otimes id) \Delta .
$$
 Then   $E_{G/N}$
and $E_{N \backslash G}$ are projections of norm one
(completely positive and completely bounded conditional expectations)
from $A_G$ onto
$ A_{N \backslash G}$ and $ A_{G/N}$ respectively
(cf. \cite{Pod95a}, as well as Proposition 2.3 and Section 6 of \cite{act}),
and
$$
{\mathcal A}_{G/N} = E_{G/N}( {\mathcal A}_G ) , \; \; \;
\text{and}
\; \; \;
{\mathcal A}_{N \backslash G}  =
E_{N \backslash G} ( {\mathcal A}_G ) .
$$
The proposition below follows immediately from the above considerations.

\bgprop
\label{density}

The *-subalgebras
$ {\mathcal A}_{N \backslash G}$ and $ {\mathcal A}_{G/N}$ are dense in 
$A_{N \backslash G}$ and $A_{G/N}$ respectively under the norm of $A_G$.

\ndprop

From \propref{density}, we have following slight reformulation of
Proposition 2.1 in \cite{simple}:

\bgprop
\label{normal-subgroup}
Let $N$ be a quantum subgroup of a compact quantum group $G$.
Then the following conditions are equivalent:

{\rm (1)}  $ A_{N \backslash G}$ is a Woronowicz $C^*$-subalgebra of $A_G$.

{\rm (1)$'$}  ${\mathcal A}_{N \backslash G}$ is a Hopf
*-subalgebra of ${\mathcal A}_G$.

{\rm (2)}  $ A_{G/N}$ is a Woronowicz $C^*$-subalgebra of $A_G$.

{\rm (2)$'$} ${\mathcal A}_{G/N}$ is a Hopf
*-subalgebra of ${\mathcal A}_G$.

{\rm (3)}  $ A_{G/N} = A_{N \backslash G}$.

{\rm (3)$'$} ${\mathcal A}_{G/N} = {\mathcal A}_{N \backslash G}$.

{\rm (4)}  $N$ is normal.
\ndprop

Because of \thmref{Hopf $*$-ideals vs. Woronowicz $C^*$-ideals},
\propref{density} and \propref{normal-subgroup}, we may (and will) work
exclusively with the dense Hopf $*$-algebras
of Woronowicz $C^*$-algebras from now on unless otherwise specified.
As remarked after the proof of Proposition 2.1 in
\cite{simple}, the counit of ${\mathcal A}_{G/N}$ is equal to the restriction morphism
$\pi|_{{\mathcal A}_{G/N}}$.

As usual, if ${\mathcal H}$ is a Hopf algebra,
${\mathcal H}^+$ denotes the
augmentation ideal (i.e. ${\mathcal H}^+$ is kernel of the counit $\varepsilon$
of ${\mathcal H}$).
Assume $N$ is a normal quantum subgroup of a compact quantum group $G$. Then
we have a Hopf $*$-algebra ${\mathcal A}_{G/{N}}$ and its
augmentation ideal ${\mathcal A}^+_{G/{N}}$.

\lmref{reconstruct-N-from-G/N} below is a key ingredient in the
proof of \thmref{thm-equivalence}, and it plays an important role in
\cite{simple} and \thmref{property-F-thm} below. In the case of
an ordinary compact group $G$, its geometric meaning
 is the trivial fact that a normal subgroup
$N$ of $G$ is the inverse image of the identity element in $G/N$ under the
quotient map. However, in the case of quantum groups using
the Hopf algebra language, it is non-trivial to prove,
especially because of related complications concerning the notion
of normality for arbitrary Hopf algebras such as
Example 1.2 in Schneider \cite{Schneid93a}.
\lmref{reconstruct-N-from-G/N} is a consequence of
Takeuchi's  Theorem 2 in \cite{Takeuchi79},
as pointed out to us by the referee,
because short exact sequences of comodules over cosemisimple Hopf algebras
are always split (cf. Theorem 3.1.5 of \cite{Dascalescu}),
from which one immediately sees using for instance
1.2.11(b) of \cite{Andrus95a} that ${\mathcal A}_G$ is faithfully
coflat over ${\mathcal A}_N$ when $N$ is normal quantum subgroup of
a compact quantum group $G$, thus the condition of
Theorem 2 in \cite{Takeuchi79} is fulfilled. We note that
cosemisimple Hopf algebras in general are not faithfully coflat
over its quotient Hopf algebras if the latter is not cosemisimple,
as Chirvasitu shows by an example in \cite{Chir11pr}. Our quotient
Hopf algebra ${\mathcal A}_N$ is, however, cosemisimple, and the
pathology in  Chirvasitu's example does not occur.

Without using of faithfully coflatness and the above references, 
an outline of another proof of \lmref{reconstruct-N-from-G/N}  
is sketched in Lemma 4.4 in \cite{simple}. Because of its usefulness and 
 for the convenience of the reader, 
we include here a detailed and self-contained proof following the 
lines in \cite{simple} (cf. 16.0.2 in Sweedler \cite{Sweedler} 
and (4.21) in Childs \cite{Childs} for finite dimensional case). 

\bglm
\label{reconstruct-N-from-G/N}
{\rm (Reconstruction of $N$ from identity in $G/N$)}

Let $(N, \pi)$ be a normal quantum subgroup of a compact quantum
group $G$. Let $\hat{\pi} = \pi|_{{\mathcal A}_G}$ be the associated morphism from
${\mathcal A}_G$ to ${\mathcal A}_N$. Then,
$$
\ker(\hat{\pi}) =
{\mathcal A}^+_{G/N} {\mathcal A}_G =
{\mathcal A}_G {\mathcal A}^+_{G/N}
= {\mathcal A}_G {\mathcal A}^+_{G/N} {\mathcal A}_G.
$$
\ndlm

\noindent
\pf
It suffices to prove
$\ker(\hat{\pi}) ={\mathcal A}^+_{G/N} {\mathcal A}_G$,
as we will have equality
$\ker(\hat{\pi}) ={\mathcal A}_G {\mathcal A}^+_{G/N}$
by the same method, and these will imply that
$$
{\mathcal A}_G {\mathcal A}^+_{G/N} {\mathcal A}_G
= {\mathcal A}^+_{G/N} {\mathcal A}_G {\mathcal A}_G
= {\mathcal A}^+_{G/N} {\mathcal A}_G
= \ker(\hat{\pi})
= {\mathcal A}_G {\mathcal A}^+_{G/N}.
$$

Consider the right ${\mathcal A}_N$-comodule structures $\alpha$ and $\beta$ on
${\mathcal A}_G$ and ${\mathcal A}_N$ defined respectively
by
$$
\alpha = (id \otimes \hat{\pi})\Delta_G:
{\mathcal A}_G  \rightarrow {\mathcal A}_G \otimes {\mathcal A}_N,
% \; \; \text{and} \; \; \;
 $$
 $$
 \beta = \Delta_N:
{\mathcal A}_N  \rightarrow {\mathcal A}_N \otimes {\mathcal A}_N ,
$$
where $\Delta_G$ and $\Delta_N$ are respectively the coproducts
of the Hopf algebras ${\mathcal A}_G$ and ${\mathcal A}_N$.
Since $\hat{\pi}$ is compatible with the coproducts,
one verifies that
$
(\hat{\pi} \otimes id) \alpha = \beta \hat{\pi}.
$
That is, the surjection $\hat{\pi}$ is a morphism of
${\mathcal A}_N$-comodules
from ${\mathcal A}_G$ to ${\mathcal A}_N$.
The Hopf algebra ${\mathcal A}_N$ is
cosemisimple by the fundamental work of Woronowicz \cite{Wor87b}
(see remarks in 2.2 of \cite{free} which assures work in
\cite{Wor87b} is valid for all compact quantum groups without separability
assumption on the underlying $C^*$-algebra $A_G$ because of
Van Daele's theorem  \cite{Daele95a} on the Haar measure based on \cite{Wor87b,Wor98a}).
Therefore it follows from Theorem 3.1.5
of \cite{Dascalescu} that every ${\mathcal A}_N$-comodule is projective.
 Hence $\hat{\pi}$ has a comodule splitting
$s:  {\mathcal A}_N \rightarrow {\mathcal A}_G $
with $\hat{\pi} s = id_{{\mathcal A}_N}$.

Let $x \in {\mathcal A}^+_{G/N}$.
It is straightforward to verify that $\pi|_{_{{\mathcal A}_{G/N}}}$ is
the counit of ${\mathcal A}_{G/N}$
(cf. remark (a) following Definition 2.2 in \cite{simple}). Hence
$\hat{\pi}(x) = 0$, and therefore
  $ {\mathcal A}^+_{G/{N}} {\mathcal A}_G \subset \ker(\hat{\pi})$.
It remains to show that
  $\ker(\hat{\pi}) \subset {\mathcal A}^+_{G/{N}} {\mathcal A}_G $.

Define a linear map $\phi$ on ${\mathcal A}_G$ by
$\phi = (s \hat{\pi}) * S = m (s \hat{\pi} \otimes S) \Delta_G$,
where $m$ and $S$ are respectively the multiplication map
and antipodal map of ${\mathcal A}_G$.
%Then using the coassociativity of $\Delta_G$ and
%$\hat{\pi} s = id_{{\mathcal A}_N}$ along with the antipodal property
%of $S$, one verifies
We show that $\phi({\mathcal A}_G) \subset {\mathcal A}_{G/N}$.
To see this, let $a \in {\mathcal A}_G$.
Using the fact that  $s$ is a comodule morphism, i.e.
$ \alpha s = (s \otimes id) \beta $ or
$(id \otimes \hat{\pi}) \Delta_G s = (s \otimes id) \Delta_N$,
along with properties of Hopf algebras morphisms, we obtain
\bgeqq
(id \otimes \pi) \Delta_G (\phi(a))
&=& (id \otimes \pi) \Delta_G ( \sum s \pi(a_{(1)}) S(a_{(2)}) ) \\
&=& (id \otimes \pi) (\sum \Delta_G (s \pi(a_{(1)}))
( S(a_{(3)}) \otimes S(a_{(2)}) ) ) \\
&=& \sum (s \otimes id) ( \Delta_N ( \pi(a_{(1)})) )
 ( S(a_{(3)}) \otimes \pi S(a_{(2)}) ) \\
&=& \sum [ s \pi(a_{(1)}) \otimes \pi(a_{(2)}) ]
 [ S(a_{(4)}) \otimes \pi S(a_{(3)}) ] \\
&=& \sum  s \pi(a_{(1)}) S(a_{(3)}) \otimes \varepsilon(a_{(2)})  \\
&=& \sum  s \pi(a_{(1)}) S(\varepsilon(a_{(2)}) a_{(3)}) \otimes 1  
= \phi(a) \otimes 1,
\ndeqq
which means that $\phi({\mathcal A}_G) \subset {\mathcal A}_{G/{N}}$, where
$\varepsilon$ is the co-unit on ${\mathcal A}_G$.

Next we observe that
$$s \hat{\pi} = (s \hat{\pi}) * \varepsilon = (s \hat{\pi}) * (S * id)
= ((s \hat{\pi}) * S ) * id = \phi * id, $$
and
\bgeqq
\varepsilon \phi(a)
&=& \varepsilon( \sum s \pi(a_{(1)}) S(a_{(2)}) )
=  \sum \varepsilon(s \pi(a_{(1)}))  \varepsilon(S(a_{(2)}) )  \\
&=&  \sum \varepsilon[ s \pi( a_{(1)} \varepsilon( a_{(2)} ) )  ]
=  \varepsilon ( s \pi(a) )  = (\varepsilon_N \pi) ( s \pi(a) )  \\
&=& \varepsilon_N ( \pi(a) ) =  \varepsilon (a),
\ndeqq
where
$\hat{\pi} s = id_{{\mathcal A}_N}$ is used.
The above means $\varepsilon \phi = \varepsilon$.
From these we obtain
\bgeqq
id-s \hat{\pi} &=& \varepsilon * id - \phi * id
= (\varepsilon - \phi) * id  \\
&=& (\varepsilon \phi - \phi) * id
= [(\varepsilon  - id) \phi] * id.
\ndeqq

Furthermore
$(id- s \hat{\pi})(a) = a$ if $a \in \ker(\hat{\pi})$, we have
$\ker(\hat{\pi}) \subset {\rm Im}(id- s \hat{\pi})$. Therefore to
show  $  \ker(\hat{\pi})  \subset {\mathcal A}^+_{G/{N}} {\mathcal A}_G$,
it suffices to show that
${\rm Im} (id- s \hat{\pi}) \subset {\mathcal A}^+_{G/{N}} {\mathcal A}_G$.
%(it will be clear that
%${\rm Im} (id- s \hat{\pi}) = {\mathcal A}^+_{G/{N}} {\mathcal A}_G$).
Since $(\varepsilon-id) \phi({\mathcal A}_G) \subset {\mathcal A}^+_{G/N}$
(because $\phi({\mathcal A}_G) \subset {\mathcal A}_{G/N}$),
the later follows from the identity
$$
id- s \hat{\pi} =
[(\varepsilon-id) \phi ] * id
=m((\varepsilon-id) \phi \otimes id)\Delta_G.
$$

This proves  \lmref{reconstruct-N-from-G/N}.
\qed

\medskip

\noindent
{\em Remarks.}
Let $\Phi$ and $\Psi$ be the notations of Schneider  \cite{Schneid93a},
$$\Phi({\mathcal A}_{G/N}):= {\mathcal A}_G {\mathcal A}_{G/N}^+ , \; \;
\text{and} \; \; \Psi(\ker(\hat{\pi})) := {\mathcal A}_{G/N}.
$$
\lmref{reconstruct-N-from-G/N}  above can be restated as saying that
for a compact quantum group $G$, the map $\Phi$ is the left inverse of $\Psi$, i.e.,
$\Phi$ is a surjection from the set of normal Hopf subalgebras of ${\mathcal A}_G$
onto the set of its normal Hopf ideals.
In addition,  Chirvasitu recently showed in \cite{Chir11pr} that
the Hopf algebra ${\mathcal A}_G$ is faithfully flat over its
Hopf subalgebras, as the author had conjectured in
an earlier version of this paper and in \cite{simple}
(cf. first part of Conjecture 1 on p3329 there).
Chirvasitu's result can be used along with those of Schneider \cite{Schneid93a}
to conclude that the map $\Phi$ is also the right inverse of $\Psi$,
i.e., $\Phi$ is also an injection,
complementing \lemref{reconstruct-N-from-G/N}.
Therefore, the maps $\Phi$ and $\Psi$ are inverses to each other.
The first result of this kind is due to Takeuchi \cite{Takeuchi72} for
{\em commutative}  Hopf algebras, using which he gave
a purely Hopf algebraic proof of the fundamental theorem
of affine algebraic group schemes \cite{DemazureGabriel70}.

In the language of Andruskiewitsch {\sl et al}
\cite{Andrus95a}, \lmref{reconstruct-N-from-G/N} implies that the sequence
\bgeqq
1 \longrightarrow
N  \longrightarrow
G \longrightarrow
G / N \longrightarrow
1,
\ndeqq
or the sequence
\bgeqq
{\mathbb C} \longrightarrow
{\mathcal A}_{G/{N}} \longrightarrow
{\mathcal A}_G \longrightarrow
{\mathcal A}_N \longrightarrow
{\mathbb C},
\ndeqq
is exact, where the one dimensional Hopf algebra ${\mathbb C}$
is the ``zero object'' in the category of Hopf algebras.
Note that in the purely algebraic
situation of Parshall and Wang \cite{ParshallWang91a},
for a given normal quantum subgroup in their sense
(i.e. a-normal as defined in our paper here),
the existence of an exact sequence is not known
and the uniqueness does not hold in general (cf. 1.6 and 6.3 loc. cit.).
\lmref{reconstruct-N-from-G/N} above shows that
such complications do not present themselves
in the world of compact quantum groups: when we have a normal quantum group,
we always have a unique exact sequence. This property might help to
formulate an appropriate notion of quantum groups in algebraic setting.

Note also that the notion of exact sequence of quantum groups in
\cite{Schneid93a} is equivalent to the notion of
strictly exact sequence in \cite{Andrus95a}
under faithful (co)flat conditions, which are fulfilled for
cosemisimple Hopf algebras and therefore for compact
quantum groups thanks to the theorem of Chirvasitu \cite{Chir11pr}
on faithfully flatness and the remarks before
 \lmref{reconstruct-N-from-G/N} on faithfully coflatness.
In general, an arbitrary Hopf algebra need not be faithfully
flat over its Hopf subalgebras, according to counter examples of
Schauenburg \cite{Schauenburg00a} constructed in response to
Question 3.5.4 of Montgomery \cite{Montgomery93a}.

We now prove the first main result of this paper.
\vv
{\em Proof of \thmref{thm-equivalence}} (cf. Schneider \cite{Schneid93a}).

(1)  $\Rightarrow$ (2) and therefore (1)  $\Rightarrow$ (3) and (1)  $\Rightarrow$ (4):
Let $(N, {\pi})$ be normal,
we show that $(N, \hat{\pi})$ is a-normal.

Let $a \in \ker(\hat{\pi})$. We show that
$ad_l(a) \in \ker(\hat{\pi}) \otimes {\mathcal A}_G$,
where
$$
ad_l(a) %:
= \sum a_{(2)} \otimes  a_{(1)} S(a_{(3)}) .
$$

By \lmref{reconstruct-N-from-G/N}
$$
\ker(\hat{\pi}) =
{\mathcal A}^+_{G/{N}} {\mathcal A}_G =
{\mathcal A}_G {\mathcal A}^+_{G/{N}}
= {\mathcal A}_G {\mathcal A}^+_{G/{N}} {\mathcal A}_G.
$$
%Set ${\mathcal I} = \ker(\hat{\pi})
%= {\mathcal A}_G {\mathcal A}^+_{G/{N}} {\mathcal A}_G.$
We assume without loss of generality $a = b c$ for
$b \in {\mathcal A}_G$ and $c \in {\mathcal A}^+_{G/{N}}$.
Then modulo $\ker(\hat{\pi}) \otimes {\mathcal A}_G$
we have
\bgeqq
ad_l(a)
&=&  \sum a_{(2)} \otimes  a_{(1)} S(a_{(3)})
= \sum b_{(2)} c_{(2)} \otimes  b_{(1)} c_{(1)} S(c_{(3)}) S(b_{(3)}) \\
&=&
\sum b_{(2)} \varepsilon(c_{(2)}) \otimes
b_{(1)} c_{(1)} S(c_{(3)}) S(b_{(3)}) \\
&=&
\sum b_{(2)} \otimes
b_{(1)} \sum ( c_{(1)} S(c_{(2)}) ) S(b_{(3)}) \\
&= &
\sum b_{(2)} \otimes
b_{(1)} \varepsilon(c) S(b_{(3)}) = 0,
\ndeqq
i.e., $ad_l(a) \in \ker(\hat{\pi}) \otimes {\mathcal A}_G$,
where the property
$$(\varepsilon- id) ({\mathcal A}_{G/N}) \subset {\mathcal A}^+_{G/N}
\subset \ker(\hat{\pi}) $$
is used, as well as the counital and antipodal properties.

Similarly, $ad_r(a) \in \ker(\hat{\pi}) \otimes {\mathcal A}_G$ by assuming
$a = cb$ for $b \in {\mathcal A}_G$ and $c \in {\mathcal A}^+_{G/{N}}$,
where
$$
 ad_r(a) %:
= \sum a_{(2)} \otimes S(a_{(1)}) a_{(3)} .
$$

Hence
$ad_l(\ker(\hat{\pi})) \subset \ker(\hat{\pi}) \otimes {\mathcal A}_G$ and
$ad_r(\ker(\hat{\pi})) \in \ker(\hat{\pi}) \otimes {\mathcal A}_G$, and $N$ is
a-normal.

(3)  $\Rightarrow$ (1): Assume $(N, \hat{\pi})$ is left a-normal.
We prove the equality ${\mathcal A}_{G/N} = {\mathcal A}_{N \backslash G}$, hence
by equivalence of (3)$'$ and (4) in \propref{normal-subgroup}, $N$ is normal.
Another proof of this is in 1.1.7 of [1], so readers familiar with it
may skip the proof below.
Note that our proof of this {\em equality} is for general Hopf algebras
with bijective antipode not necessarily associated with
compact quantum groups, just as  1.1.7 of [1].

Let $a \in {\mathcal A}_{N \backslash G}$, i.e.,
$$
\sum a_{(1)} \otimes a_{(2)} - 1 \otimes a
\in \ker(\hat{\pi}) \otimes {\mathcal A}_G.
$$
 Applying $(ad_l \otimes id)$ to this and using the condition in
the definition of left a-normal, we obtain
$$
\sum a_{(2)} \otimes  a_{(1)} S(a_{(3)}) \otimes a_{(4)}
- 1 \otimes 1 \otimes a \in
 \ker(\hat{\pi}) \otimes {\mathcal A}_G \otimes {\mathcal A}_G.
$$
Multiplying the second and the third factors,
%i.e., applying $id \otimes m$ to both sides of the last expression,
we obtain
$$
\sum a_{(2)} \otimes  a_{(1)} S(a_{(3)}) a_{(4)} - 1 \otimes a \in
 \ker(\hat{\pi}) \otimes {\mathcal A}_G .
$$
By the antipodal and counital properties,
$$
\sum a_{(2)} \otimes  a_{(1)} S(a_{(3)}) a_{(4)}
=  \sum a_{(2)} \otimes  a_{(1)} \varepsilon(a_{(3)})
=  \sum a_{(2)} \otimes  a_{(1)} .
$$
Hence
$$
 \sum a_{(2)} \otimes  a_{(1)}
- 1 \otimes a \in
 \ker(\hat{\pi}) \otimes {\mathcal A}_G ,
\; \; \text{and therefore}$$
$$
 \sum a_{(1)} \otimes  a_{(2)}
- a \otimes 1 \in
  {\mathcal A}_G \otimes \ker(\hat{\pi}), \; \;
\text{i.e.,} \; \;
 a \in {\mathcal A}_{G/N}.
$$

That is, we have an inclusion ${\mathcal A}_{N \backslash G} \subset  {\mathcal A}_{G/N}$.
Using definition of
${N \backslash G}$ and $G/N$ and properties of the antipode $S$
we immediately have that
$S({\mathcal A}_{N \backslash G}) = {\mathcal A}_{G/N}$ and
$S({\mathcal A}_{G/N}) = {\mathcal A}_{N \backslash G}$.
Using this and applying $S$ to the above inclusion
we obtain
${\mathcal A}_{G/N} = S({\mathcal A}_{N \backslash G})
\subset S({\mathcal A}_{G/N}) = {\mathcal A}_{N \backslash G},$
i.e., we also have ${\mathcal A}_{G/N} \subset {\mathcal A}_{N \backslash G}$.
Hence ${\mathcal A}_{G/N} = {\mathcal A}_{N \backslash G}$.

(4) $\Rightarrow$ (1): The proof is similar to
(3)  $\Rightarrow$ (1) above.

(2) $\Rightarrow$ (1): This follows from either (4) $\Rightarrow$ (1) or
(3)  $\Rightarrow$ (1).
\qed
\vv
{\em Remarks.}
Although for general Hopf algebra not necessarily associated with compact quantum groups,
${\mathcal A}_{G/N} = {\mathcal A}_{N \backslash G}$ follows from
{\em either} $(N, \hat{\pi})$ being left a-normal or
right a-normal, if the assertion in \lmref{reconstruct-N-from-G/N}
is not valid for such Hopf algebra, which is a key ingredient
in the proof of the implication (1) $\Rightarrow$ (2) above,
we probably cannot expect left a-normal or right a-normal to follow from
${\mathcal A}_{G/N} = {\mathcal A}_{N \backslash G}$.
We note that no example seems to be known of an
algebraic quantum group that is left a-normal
but not right a-normal, or vice versa.
We suspect such an example may come from a Hopf algebra
not faithfully coflat over its quotient Hopf algebra.
We are not aware if an example of the latter
has been produced in Hopf algebra literature,
which should exist in view of counter examples
on faithfully flatness (cf. \cite{Schauenburg00a}).

\medspace
\medspace

The Equivalence \thmref{thm-equivalence}
enables results on normality in Hopf algebra literature
be applicable to normal quantum subgroups of compact quantum groups
in our sense, such as those in
\cite{ParshallWang91a, Schneid93a, Takeuchi94a,Andrus95a},
noting that conormal in \cite{Schneid93a} and
a-normal are the same concept.

\section{Third Fundamental Isomorphism Theorem for Quantum Groups
and Properties $F$ and $FD$}
\label{applications}

In this section, we give several applications of \thmref{thm-equivalence}.
Because of this theorem and other results in Sections 2 and 3,
we mostly focus on the dense Hopf $*$-subalgebras
associated to compact quantum groups in this section.
For ease of notation, we now use undecorated $\pi$ for Hopf algebra
morphism ${\mathcal A}_G \rightarrow {\mathcal A}_N$,
omitting the hat in $\hat{\pi}$ whenever no confusion arises.

The three fundamental isomorphism theorems in the
theory of groups are foundational results on structure of groups.
One way naturally expect their analogs to be
valid in the theory of quantum groups. Unfortunately, quantum analog of
the first fundamental isomorphism theorem is not always true for epimorphism
of quantum groups (i.e. injection of Hopf algebras) except for the situation where
exact sequence can be constructed,
cf. \cite{ParshallWang91a,Schneid93a,Takeuchi94a,Andrus95a}.
For instance, not every Woronowicz
$C^*$-subalgebra of $A_G$ is of the form $A_{G/N}$ with $N$ a normal quantum
subgroup of $G$. This already fails when $A_G$ is the group $C^*$-algebra $C^*(F_2)$
of the free group $F_2$ on two generators, because a Woronowicz $C^*$-subalgebra
of $A_G$ is not of the form  $A_{G/N}$ unless
it is the group $C^*$-algebra $C^*(\Gamma)$ of a normal
subgroup $\Gamma$ of $F_2$. Finally, it is not clear how a quantum analog of
the second fundamental isomorphism theorem can be formulated.

However, on the bright side, as an application of \thmref{thm-equivalence},
we have the following complete analog of the {\em Third Fundamental Isomorphism Theorem}
for compact quantum groups.

\bgth
{\rm (Third Fundamental Isomorphism Theorem)}
\label{3rd-fundamental-isomorphism-thm}
Let $(N, \pi)$ be a normal quantum subgroup of $G$. Let $(H, \theta)$ be a
quantum subgroup of $G$ that contains $(N, \pi)$, i.e., there is a morphism
$\pi_1$ from ${\mathcal A}_{H}$ to ${\mathcal A}_{N}$ such that
$(N, \pi_1)$ is a quantum subgroup of $H$
with $\pi = \pi_1 \theta$. Then $(N, \pi_1)$ is normal in $H$.
If furthermore $H$ is normal in $G$
and letting $\theta' = \theta |_{{\mathcal A}_{G/N}}$,
the restriction of $\theta$ to ${\mathcal A}_{G/N}$,
then $(H/N, \theta')$ is normal in $G/N$ and
$$
{\mathcal A}_{(G/N)/(H/N)} = {\mathcal A}_{ G/H}
= A_{ H \backslash G} = A_{(N \backslash H)  \backslash (N \backslash G)}.
$$
\ndth
\noindent
\pf
Let $z \in {\mathcal A}_H$ be such that $\pi_1(z)=0$.
Assume $z = \theta(x)$
for some $x \in {\mathcal A}_G$. Then $\pi(x) =0$.
Since $N$ is a-normal in $G$ by \thmref{thm-equivalence},
we have
$$
\sum \pi(x_{(2)}) \otimes S_G(x_{(1)}) x_{(3)} = 0,
\; \;
\sum \pi(x_{(2)}) \otimes x_{(1)} S_G(x_{(3)}) = 0.
$$
Hence
$$
\sum \pi_1 \theta(x_{(2)}) \otimes \theta (S_G(x_{(1)}) x_{(3)}) = 0,
\; \;
\sum \pi_1 \theta(x_{(2)}) \otimes \theta (x_{(1)} S_G(x_{(3)})) = 0.
$$
It follows that
$$
\sum \pi_1 (z_{(2)}) \otimes S_H(z_{(1)}) z_{(3)} = 0,
\; \;
\sum \pi_1 (z_{(2)}) \otimes z_{(1)} S_H(z_{(3)}) = 0.
$$
This means that $(\pi_1, N)$ is a-normal and is therefore normal
in $H$ by \thmref{thm-equivalence}.

Let $a \in {\mathcal A}_{ G/N}$.  Then it is immediate to verify that
$\theta(a) \in {\mathcal A}_{ H/N}$, and therefore
$\theta({\mathcal A}_{ G/N})$ is contained in $ {\mathcal A}_{ H/N}$.
Conversely, let $b \in {\mathcal A}_{ H/N}$. Assume $b = \theta(a)$
for some $a \in {\mathcal A}_{G}$.
Put
$$\bar{a} = E_{G/N}(a) = (id \otimes h_N \pi) \Delta_G(a). $$
Then $\bar{a} \in {\mathcal A}_{ G/N}$ and
 $$b = E_{H/N}(b) = (id \otimes h_N \pi_1) \Delta_H(b)$$
 by the remarks before \propref{density} applied to $G/N$ and $H/N$ respectively.
 Moreover, we have
\bgeqq
\theta(\bar{a}) &=& (\theta  \otimes h_N \pi) \Delta_G(a)
= (\theta  \otimes h_N \pi_1 \theta) \Delta_G(a)  \\
&=& (id  \otimes h_N \pi) \Delta_H(\theta(a)) =(id \otimes h_N \pi_1) \Delta_H(b) \\
&=& b.
\ndeqq
That is
$\theta(\bar{a}) = b$. Therefore $\theta({\mathcal A}_{ G/N})= {\mathcal A}_{ H/N}$.

For ease of notation, let $G' = G/N$
and $H' = H/N$. The above shows that $(H', \theta')$ is a quantum subgroup
of $G'$.

Now assume $(H, \theta)$ is normal. If $a \in {\mathcal A}_{G'/H'}$, that is,
$(id \otimes \theta') \Delta(a) = a \otimes 1_{H'}$, then
it is clear that $a $ is in ${\mathcal A}_{G/H}$ since
$\theta' = \theta |_{A_{G/N}}$ and $1_{H'}=1_H$.
Conversely, if $a \in {\mathcal A}_{G/H}$, that is
$(id \otimes \theta) \Delta(a) = a \otimes 1_{H}$,
then
$(id \otimes \pi_1 \theta) \Delta(a) = a \otimes 1_{N}$.
This means that
$a \in {\mathcal A}_{G'} =  {\mathcal A}_{G/N}$.
Since the coproduct for the quantum group ${G'} = {G/N}$ is a restriction of
$\Delta$,  we have
$$(id \otimes \theta') \Delta(a) =
(id \otimes \theta) \Delta(a) = a \otimes 1_{H} =  a \otimes 1_{H'}.$$
Hence  $a \in {\mathcal A}_{G'/H'}$ and
$A_{G'/H'} = A_{ G/H}.$

The result is completely proved.
\qed
\vv
{\em Remark:}
Instead of an isomorphism such as in
${(G/N)/(H/N)} \cong { G/H}$ in group theory,
we have {\em exact equalities} of quantum function algebras
in \thmref{3rd-fundamental-isomorphism-thm} above.
\vv
The proof of \thmref{3rd-fundamental-isomorphism-thm}
actually yields the following stronger result without assuming
$(H, \theta)$ to be normal, which should be useful in
harmonic analysis on homogeneous spaces.
\vv
{\bf \thmref{3rd-fundamental-isomorphism-thm}$'$.} {\em
Let $(N, \pi)$ be a  normal quantum subgroup of $G$. Let $(H, \theta)$ be a
(not necessarily normal) quantum subgroup of $G$  that contains $(N, \pi)$,
i.e., there is a morphism $\pi_1$ such that $(N, \pi_1)$ is
a quantum subgroup of $H$ with $\pi = \pi_1 \theta$.
Let $\theta' = \theta |_{{\mathcal A}_{G/N}}$,
the restriction of $\theta$ to ${\mathcal A}_{G/N}$.
Then $(N, \pi_1)$ is normal in $H$, $(H/N, \theta')$ is a
quantum subgroup of $G/N$, and
$$
A_{(G/N)/(H/N)} = A_{ G/H},
\; \; \;
\; \; \;
A_{(N \backslash H)  \backslash (N \backslash G)}
= A_{ H \backslash G}.
$$
}
\indent
For other applications of \thmref{thm-equivalence}, we
first recall the following properties of
compact quantum groups  (cf. \cite{simple}).

\bgdf
\label{property-F-etc}
A compact quantum group $G$ is said to have {\bf property $F$} if each
Woronowicz $C^*$-subalgebra of $A_G$ is of the form $A_{G/N}$ for some
normal quantum subgroup $N$ of $G$;
$G$ is said to have {\bf property $FD$}
if each quantum subgroup of $G$ is normal.
\nddf

These notions are motivated by the following facts
(see Propositions 2.3 and 2.4 in \cite{simple}): if $G$ is a
compact group, then its function algebra $C(G)=A_G$ has property $F$;
if $G$ is the dual of a discrete group $\Gamma$, then
its quantum function algebra $C^*(\Gamma)$ has property $FD$.
Therefore compact quantum groups with property $F$
are closest to compact groups,
while compact quantum groups with property $FD$ are
closest to the compact quantum group dual
of discrete groups.

We note that the notions {property $F$} and {property $FD$} above
can be defined almost verbatim for all Hopf algebras -- one only needs to replace
the words ``compact quantum group'' (resp. ``Woronowicz $C^*$-subalgebra'')
with the words ``Hopf algebra'' (resp. ``Hopf subalgebra'') in
the above definition. Moreover,
because of \thmref{thm-equivalence} and remarks following
\lmref{reconstruct-N-from-G/N}, in Hopf algebra language,
$G$ has {property $F$} if each Hopf subalgebra of ${\mathcal A}_G$ is normal
in the sense of 3.4.1 in \cite{Montgomery93a};
it has {$FD$} if each Hopf $*$-ideal of ${\mathcal A}_G$ is normal
in the sense of 3.4.5 in \cite{Montgomery93a}. It follows from
discussions after 3.4.5 in \cite{Montgomery93a} that if
$G$ is a finite quantum group, i.e. ${\mathcal A}_G$ is finite dimensional,
$G$ has property $F$ if and only if its Pontryagin dual $G_d$
has property $FD$, where the Pontryagin dual $G_d$
is the finite quantum group with quantum function algebra equal to
the dual Hopf algebra ${\mathcal A}^\prime_G$ of ${\mathcal A}_G$.

In Franz {\sl et al.} \cite{FrSkTo12}, in different terminology, 
compact quantum groups with property $FD$ are called {\bf hamiltonian}, of which 
quantum groups in the $DS$ family are a special case. 
An example of noncommutative and noncocommutative quantum group in 
the $DS$ family (and therefore an example with property $FD$) 
is given in section 6 of \cite{FrSkTo12}. 

In our earlier work \cite{simple}, it is shown that
all the quantum groups obtained by deformation of
compact Lie groups, such as the compact real forms of Drinfeld-Jimbo
quantum group and Rieffel's deformation of compact Lie groups,
and all universal quantum groups \cite{free,uqg,qsym} (except $A_u(Q)$)
have {property $F$}. These are non-trivial and natural
examples of compact quantum groups that have property $F$.
Despite this multitude of examples, as discussed in
\cite{simple,simpleprob}, much more is to be explored
and it would be important to develop a classification theory
for simple quantum groups with property $F$.

Furthermore, we have following result on the
structure of compact quantum group with property $F$.

\bgth
\label{property-F-thm}
Quotient group $G/N$ by a normal quantum subgroup
$(N, \pi)$ has property $F$ if 
$G$ is a compact quantum group with property $F$.
\ndth

\noindent
\pf
Let $(N, \pi)$ be a normal
quantum subgroup of $G$ and let
${\mathcal C} \subset {\mathcal A}_{G/N}$ be a
Hopf subalgebra of ${\mathcal A}_{G/N}$.
We show that there is a normal quantum subgroup $(K, \theta')$
of $G/N$ such that ${\mathcal C} = {\mathcal A}_{(G/N)/K}$.
The proof below is suggested by the referee, replacing our original
longer and direct proof without using
\thmref{3rd-fundamental-isomorphism-thm} and \lmref{reconstruct-N-from-G/N}.

As ${\mathcal C} \subset {\mathcal A}_{G/N}$ is also
a Hopf subalgebra of ${\mathcal A}_G$, by property $F$ of $G$,
let $(H, \theta)$ be a
normal quantum subgroup of $G$ such that
${\mathcal C} = {\mathcal A}_{G/H}$.
By \thmref{3rd-fundamental-isomorphism-thm}, it is enough
to show that $(H, \theta)$ contains $(N, \pi)$, i.e.
$\ker(\theta) \subset \ker({\pi})$, because then
the normal quantum subgroup $(K, \theta')$
of $G/N$ we need is simply $(H/N, \theta')$
in \thmref{3rd-fundamental-isomorphism-thm}:
${\mathcal C} = {\mathcal A}_{G/H} = {\mathcal A}_{(G/N)/(H/N)}. $

By \lmref{reconstruct-N-from-G/N},
(noting that we have omitted the hats to simplify notation)
we have
$$\ker(\theta) =  {\mathcal A}_G {\mathcal A}^+_{G/H} {\mathcal A}_G,
\; \;
\ker({\pi}) = {\mathcal A}_G {\mathcal A}^+_{G/N} {\mathcal A}_G.
$$
Since ${\mathcal A}_{G/H} \subset  {\mathcal A}_{G/N}$,  we have
$\ker(\theta) \subset \ker({\pi})$.
\qed
\vv
{\em Remarks.}

(1) 
By considering formal dual to \thmref{property-FD-thm}.(2) below, 
a quantum subgroup $(H, \theta)$ 
of a compact quantum group $G$ with property $F$ probably does not  
always have property $F$, though it does if it is assumed that 
inverse images of Hopf subalgebras of ${\mathcal A}_H$ under
$\theta$ are Hopf subalgebras of  ${\mathcal A}_G$. 
This assumption turns out to be trivial in the sense that 
 $(H, \theta)$ is either the identity group or the full group $G$, 
 as kindly pointed out to us by the referee. 
Relevant parts of Theorem 3.7 in \cite{simple} and 
Theorem 5.6 in \cite{simpleprob} on the same result 
need to be modified in the form stated in \thmref{property-F-thm} above.

(2)
We note that though there are more Woronowicz $C^*$-ideals
in $A_G$ than Hopf $*$-ideals in ${\mathcal A}_G$ in the correspondence of
\thmref{Hopf $*$-ideals vs. Woronowicz $C^*$-ideals}
(see Remark (a) after Lemma 4.3 in \cite{simple}),
Woronowicz $C^*$-subalgebras of $A_G$ and Hopf
subalgebras\footnote{Note that Woronowicz's fundamental work \cite{Wor87b,Wor98a} implies that
a Hopf subalgebra of ${\mathcal A}_G$ is automatically closed under the *-operation since
it is invariant under the antipode.}
of ${\mathcal A}_G$ are in bijective correspondence:
every Woronowicz $C^*$-subalgebra $B$ of $A_G$ uniquely corresponds to
its canonical dense Hopf $*$-subalgebra $\mathcal B$ of $B$.
In connection with this, it would be of interest to answer the following question,
which seems to have an affirmative answer according to our preliminary
investigation for several special cases.

\bgques
\label{subalg}
Is a Woronowicz $C^*$-subalgebras of a full Woronowicz
$C^*$-algebra necessarily full?
\ndques

A related question is the following one on the relation between a
Woronowicz $C^*$-ideal of a Woronowicz $C^*$-subalgebra and the
ideal it generates in the original Woronowicz $C^*$-algebra.
We formulate two equivalent versions:

\bgques
\label{pullback ideal}
Assume $A$ is a Woronowicz $C^*$-algebra and $A_0$ a Woronowicz $C^*$-subalgebra with  Hopf $*$-subalgebras ${\mathcal A}$ and ${\mathcal A}_0$ respectively.

{
{\rm (a)} Let $I_0$ be a Woronowicz $C^*$-ideal of $A_0$. Let $I : = A I_0 A$ be the Woronowicz
$C^*$-ideal of $A$ generated by $I_0$. Is the identity $I_0 = I \cap  A_0$ always true?

{\rm (b)} Let ${\mathcal I}_0$ be a Hopf $*$-ideal of ${\mathcal A}_0$.
Let ${\mathcal I} : = {\mathcal A} {\mathcal I}_0 {\mathcal A}$ be the Hopf $*$-ideal
of ${\mathcal A}$ generated by ${\mathcal I}_0$.
Is the identity ${\mathcal I}_0 = {\mathcal I} \cap {\mathcal A}_0$ always true?
}
\ndques

Parts (a) and (b) of \quesref{pullback ideal} are equivalent by
\thmref{Hopf $*$-ideals vs. Woronowicz $C^*$-ideals}.

It turns out the answer to the above question is affirmative for
some quantum groups but negative for others. 
Before turning to examples and counter examples kindly
provided to us by the referee, we recall the following
closely related fact pointed out to us also by the referee.

\bgprop
\label{faithfully flat algebras}
Let ${\mathcal A}_0$ be a subalgebra of an algebra ${\mathcal A}$ over the
complex numbers. Then the following are equivalent.

 {\rm (1)} ${\mathcal A}$  is left faithfully flat over ${\mathcal A}_0$.

 {\rm (2)} ${\mathcal A}$  is left flat over ${\mathcal A}_0$ and for
 every left ideal ${\mathcal I}_0$ in ${\mathcal A}_0$,
  ${\mathcal I}_0 =  {\mathcal A} {\mathcal I_0} \cap {\mathcal A}_0$.
\ndprop

See p33 of Bourbaki \cite{Bourbaki72} for a proof.
It is clear that a similar result relating right faithfully flatness and
right ideal is also valid.

\medspace
\medspace
\noindent
{\em Example.}
By Chirvasitu's theorem \cite{Chir11pr}, a cosemisimple Hopf algebra is
both left and right faithfully flat over its Hopf subalgebras.
It follows then from \propref{faithfully flat algebras} that
\quesref{pullback ideal} has an affirmative answer for
compact groups because the relevant algebras are commutative and left ideals
are two sided ideals. In the notation of the question above, here $A=C(G)$,
the continuous function algebra on a compact group $G$,
$A_0$ is the algebra of functions on a quotient group $G/N$ by a
normal subgroup, $I_0$ defines a subgroup $H/N$ of the quotient group,
and $I$ defines the pullback $H$ in $G$ of $H/N$ under the quotient map
from $G$ to $G/N$ .

\medspace
In view of the discussions above, the following definition is natural.

\bgdf
\label{pullback  property}
We say
\iffalse a Woronowicz $C^*$-algebra \fi
a compact quantum group $G$
has the {\bf pullback property} if
the answer for \quesref{pullback ideal} is affirmative for $A = A_G$,
equivalently for ${\mathcal A} = {\mathcal A}_G$.
\iffalse
for every normal quantum subgroup $N$
of $G$ and every quantum subgroup $K$ of $G/N$,
the answer for \quesref{pullback ideal} is affirmative.
and the morphism $\rho$ is a surjection from ${\mathcal A}_{K}$ onto ${\mathcal A}_{H/N}$.
\fi
\nddf

Note the above {\bf pullback property} can be defined for
inclusion of rings and ($C^*$-)algebras not necessarily associated with quantum groups
when the words ``Hopf $*$-ideal'' or ``Woronowicz $C^*$-ideal'' in
\quesref{pullback ideal} are substituted by ``ideal''.
The following example shows that, unlike function algebra over compact groups
in the example above, group ($C^*$-)algebras do not have pullback property
in general.

\medspace
\medspace
\noindent
{\em Counter Example.}
Let $A = {\mathcal A} = {\mathbb C} S_3$ be the group ($C^*$-)algebra
of the symmetric group $S_3$ on three symbols, and
$A_0 = {\mathcal A}_0 = {\mathbb C} \Gamma \cong {\mathbb C} {\mathbb Z} /3{\mathbb Z}$
the group  ($C^*$-)algebra of the alternating subgroup $\Gamma$ of $S_3$
generated by the three cycle $(123)$.
We claim that there exists an ideal $I_0$ in $A_0$ for which the answer
to \quesref{pullback ideal} is negative.

To see this, using (non-commutative) Fourier transforms, identify
 $A$ with $ {\mathbb C} \oplus {\mathbb C} \oplus B$,
 and $A_0$ with ${\mathbb C} \oplus B_0$,
 where $B = M_2({\mathbb C})$ is the $2 \times 2$ matrix algebra,
 corresponding to the (unique) two dimensional irreducible representation of $S_3$,
 and $B_0= {\mathbb C} \oplus {\mathbb C}$, corresponding to the
two non-trivial irreducible one dimensional representations of $\Gamma$.
Since the restriction of the two dimensional irreducible representation
of $S_3$ to $\Gamma$ is a direct sum of the two non-trivial irreducible
one dimensional representations of $\Gamma$ (see e.g. \cite{Mackey76}, p150),
we see that under the above Fourier transforms,
$B_0$ is included as the subalgebra of diagonal
matrices in $B$ (note that $B$ is not a Hopf algebra).
The inclusion  $B_0 \subset B$ does not have pullback property
because $B$ has no nonzero proper ideal
and for the nonzero proper ideal $I_0 : = {\mathbb C} \oplus 0$ in $B_0$,
$ BI_0B \cap B_0 =  B_0  \neq I_0 $. Since $I_0$ is also an
ideal in $A_0$, we have $ A I_0 A \cap A_0 = BI_0B \cap B_0 =  B_0  \neq I_0 $.

\medspace
\medspace
\noindent
We are ready for the following result on property $FD$.
\bgth
\label{property-FD-thm}
Let $G$ be a compact quantum group with property $FD$.
Then

{\rm (1)} quantum subgroup $(H, \theta)$ of $G$  also has property $FD$;

{\rm (2)} quotient group $G/N$
by a normal quantum subgroup $(N, \pi)$ also has property $FD$ provided $G$
has pullback property.
\ndth
\noindent
\pf
(1) Let $(H_1, \pi_1)$ be a quantum subgroup of $H$.
Then $(H_1, \pi_1 \theta)$ is a quantum subgroup of $G$
that is contained in $(H, \theta)$.
By property $FD$, $(H_1, \pi_1 \theta)$ is normal in $G$.
By \thmref{3rd-fundamental-isomorphism-thm},
$(H_1, \pi_1)$ is normal in $H$.
This shows that $H$ has property $FD$.

(2)
Assume $G$ is a compact quantum group with the
pullback property in additional to property $FD$.

Let $(N , \pi)$ be a normal quantum subgroup of $G$ and
let $(K, \pi_0)$ be a quantum subgroup of $G/N$.
Let ${\mathcal I}_0$ be the kernel of ${\pi}_0$ in ${\mathcal A}_{G/N}$ and
${\mathcal I}_N$ the kernel of ${\pi}$ in ${\mathcal A}_{G}$.
Identify ${\mathcal A}_{K}$  with $({\mathcal A}_{G/N})/{\mathcal I}_0$.
Put ${\mathcal I} : = {\mathcal A} {\mathcal I}_0 {\mathcal A}$.
\iffalse $I = A_G I_0 A_G$. \fi
Then ${\mathcal I}$ is a Hopf ($*$-)ideal in ${\mathcal A}_G$ and defines a
quantum subgroup $(H, \theta)$ of $G$, where
$${\theta}: {\mathcal A}_G \longrightarrow {\mathcal A}_G / {\mathcal I},
\; \; \;
{\mathcal A}_H :=  {\mathcal A}_G / {\mathcal I}. $$

By the co-unital property $\varepsilon_{G/N} = \varepsilon_{K} \pi_0$
and the fact that $\varepsilon_{G/N}$ is the restriction of $\pi$
to ${\mathcal A}_{G/N}$, we have $\pi = \varepsilon_{K} \pi_0$ on
 ${\mathcal A}_{G/N}$. It follows that
 ${\mathcal I}_0 \subset {\mathcal I}_N$ and therefore we have
 an inclusion ${\mathcal I} \subset {\mathcal I}_N$.
 This means that $(H, \theta)$ is a quantum subgroup of
 $G$ containing $(N, \pi)$ and $(N, \pi_1)$ is normal in $H$ as shown
 in \thmref{3rd-fundamental-isomorphism-thm}, where $\pi_1$ is defined
 by $\pi_1(\theta(a)) : = \pi(a)$, for $a \in {\mathcal A}_{G}$.
 On the other hand,
 every element in ${\mathcal A}_{K}$ is of the form $\pi_0(a)$
for some $a \in {\mathcal A}_{G/N}$. If $\pi_0(a)=0$, then
$\theta(a) =0$ because ${\mathcal I}_0$ is contained in
${\mathcal I} = {\mathcal A} {\mathcal I}_0 {\mathcal A}$.
Hence
$$\rho(\pi_0(a)) : = \theta(a)$$
gives a well defined morphism from
${\mathcal A}_{K}$ to ${\mathcal A}_H $, where $\pi_0(a) \in {\mathcal A}_{K}$.
It can be checked that $\rho$ is  a morphism of Hopf algebras.

We summarize all the morphisms in the following commutative diagram,
where horizontal sequences are exact in the sense of \cite{Andrus95a}.
\medspace
\medspace
\bgeqq
\begin{CD}
{\mathbb C} @>>> {\mathcal{A}_{G/N}}  @>>> {\mathcal{A}_G}
  @> {\pi}>>  {\mathcal{A}_N}  @>\varepsilon_N>> {\mathbb C}  \\
  @.@VV{{\pi}_0}V@VV{{\theta}}V@|@.   \\
{\mathbb C}  @>>> {\mathcal{A}_{K}}  @>{\rho}>> {\mathcal{A}_H}
  @>{{\pi}_1}>>  {\mathcal{A}_N}  @>\varepsilon_N>> {\mathbb C}
\end{CD}
\ndeqq
\medspace
\medspace

Moreover, the image of $\rho$ is in ${\mathcal A}_{H/N}$ because by
property of Hopf algebra morphisms, for $a \in {\mathcal A}_{G/N}$ we have
\bgeqq
& & (id_H \otimes \pi_1) \Delta_H (\rho(\pi_0(a)))
= (id_H \otimes \pi_1) \Delta_H (\theta(a)) \\
&=&  (\theta \otimes \pi_1 \theta) \Delta_G (a)
= (\theta \otimes id_N) ( id_G \otimes \pi) \Delta_G (a)  \\
&=&  (\theta \otimes id_N)  (a \otimes 1_N)
=  \theta(a) \otimes 1_N = \rho(\pi_0(a)) \otimes 1_N,
\ndeqq
which means $\rho(\pi_0(a)) \in {\mathcal A}_{H/N}$.

In addition, if $\theta(a) = 0$ for some $a$ in ${\mathcal A}_{G/N}$,
then $a$ is ${\mathcal I} \cap {\mathcal A}_{G/N}$,
which is  ${\mathcal I}_0$ by assumption, and $\pi_0(a) = 0$.
Hence $\ker(\rho) = 0$ and $\rho$ is an injection.
\iffalse
Conversely, let $b \in {\mathcal A}_{ H/N}$. Then the proof
in \thmref{3rd-fundamental-isomorphism-thm} shows that
$b = \theta(\bar{a})$ if  $b = \theta(a)$ for some $a$ in
$\mathcal{A}_G$ and $\bar{a} = E_{G/N}(a) = (id \otimes h_N \pi) \Delta_G(a)$.
\fi
As in the proof of
\thmref{3rd-fundamental-isomorphism-thm},
$\theta({\mathcal A}_{ G/N})= {\mathcal A}_{ H/N}$.
Therefore, $\rho$ is also a surjection onto  ${\mathcal A}_{H/N}$.
Hence we have an identity (and isomorphism)
$\rho({\mathcal A}_{K}) = {\mathcal A}_{H/N}$.

Since $H$  is normal in $G$ by property $FD$ of $G$,
\thmref{3rd-fundamental-isomorphism-thm} guarantees that
${K}$ (which is ${H/N}$) is normal in $G/N$ and the proof is complete.
\qed

\medspace
\medspace
\noindent
{\em Remark.} An examination of the proofs of \thmref{thm-equivalence}
and results in this and the next section shows that they are also valid for
cosemisimple Hopf algebras when the statements are appropriately modified.

\section{Other Properties of Normal Quantum Subgroups}
\label{properties-of-normal}

Tensor products of $C^*$-algebras is the analog of product of locally
compact spaces. Unlike the classical situation of spaces,
the algebraic tensor product $A_1 \otimes_{alg} A_2$ of two $C^*$-algebras $A_1$ and $A_2$
may have more than one $C^*$-norm. It has two canonical $C^*$-norms that may not
agree, the {\em maximal} one and the {\em minimal} one. Its
$C^*$-algebraic completions under these two norms are denoted respectively by
$A_1 \otimes_{max} A_2$ and $A_1 \otimes_{min} A_2$.
\iffalse
 The minimal $C^*$-tensor product $A_1 \otimes_{min} A_2$ is also called
 the {\em spatial} tensor product and is denoted simply by
 $A_1 \otimes A_2$, because it coincides with
 the completion of the algebraic tensor product $A_1 \otimes_{alg} A_2$
 under the operator norm when  $A_1 \otimes_{alg} A_2$ is faithfully
 represented on a tensor product of two Hilbert spaces.
 \fi
The maximal tensor product $A_1 \otimes_{max} A_2$ of two
full Woronowicz $C^*$-algebras $A_1 = A_{G_1}$ and $A_2 = A_{G_2}$
is also full \cite{prod}. We use $G_1 \times G_2$ to denote
the corresponding compact quantum group.

\bgprop
\label{product-group}
Let $G_1, G_2$ be compact quantum groups with normal quantum subgroups
$(N_1, \pi_1)$ and $(N_2, \pi_2)$.
Then $N: = (N_1 \times N_2, \pi_1 \otimes \pi_2)$ is a normal quantum subgroup of
$G:= G_1 \times G_2$  and $G/N \cong G_1/{N_1} \times G_2/{N_2} $,
where $\pi_1 \otimes \pi_2$ denotes the corresponding tensor product morphism
$$
\pi_1 \otimes \pi_2: A_{G_1} \otimes_{max} A_{G_2} \longrightarrow
A_{N_1} \otimes_{max} A_{N_2}.
$$
\ndprop

\noindent
\pf
According to \secref{notions-of-normal} and
 \secref{proof-equivalence}, and using the formula
 of the Haar measure on $G_1 \times G_2$ and
 formula for the coproduct of the tensor product in \cite{prod},
 one obtains immediately
$$\displaystyle{
{\mathcal A}_{G/{N}} = E_{G/{N}} ( {\mathcal A}_{G_1} \otimes {\mathcal A}_{G_2} )=
 {\mathcal A}_{G_1/{N_1}} \otimes {\mathcal A}_{G_2/{N_2}}
              } .
$$
This is a Hopf subalgebra of
${\mathcal A}_{G_1} \otimes {\mathcal A}_{G_2} $. Hence
the proposition follows from the equivalence of (3)$^\prime$
and (4) in \propref{normal-subgroup}.
\qed

\medskip

An appropriate reformulation of \propref{product-group} for minimal tensor product is
also valid when one of  the $C^*$-algebras $A_{G_1}$, $A_{G_2}$, $A_{N_1}$ and  $A_{N_2}$
is exact. We will not elaborate on this except recalling that
a $C^*$-algebra $A$ is called {\bf exact} if the functor $A \otimes_{min} ? $
preserves short exact sequences.

\medskip

For free product \cite{free}, the situation is quite different from tensor product.
The naive analog of \propref{product-group} is false, even for the
simple case with $\pi_1 = id_1$ and $\pi_2 = \varepsilon_2$, as described
in the following proposition.

\bgprop
\label{normalfreeproduct}
Let $G_1, G_2$ be compact quantum groups (not necessarily duals of discrete groups).
Let $G= G_1 \hat{*} G_2$ be the  free product compact quantum group
underlying
$A_{G_1} * A_{G_2}$.
Let $\pi$
be the natural embedding of $G_1$ into $G$
defined by the surjection
$$\pi: A_{G_1} * A_{G_2} \longrightarrow A_{G_1}, $$
$$\pi = id_1 * \varepsilon_2 .  $$
If $G_1$ has at least one irreducible representation
of dimension greater than one,
then $(G_1, \pi)$ is {\bf not} a normal quantum subgroup of
$G_1 \hat{*} G_2$.
Otherwise,  $(G_1, \pi)$ is
normal in $G_1 \hat{*} G_2$.
\ndprop
The hat in the symbol $\hat{*}$ above
signifies the ``Fourier transform'' of $*$.
\iffalse
The above result is obtained from
\propref{product-group} by setting with $\pi_1 = id_1$ and
$\pi_2 = \varepsilon_2$, but applied to free product.
\fi

\medspace
\medspace
\noindent
\pf
Let
$$u = \sum_{ij} e^u_{ij} \otimes u_{ij}, \; \; \;
v = \sum_{kl} e^v_{kl} \otimes v_{kl}
$$
be irreducible representations of $G_1,$ $ G_2$ respectively.
Assume that the dimension of $u$ is greater than one.
Let $\bar{u} = \sum_{ij} e^u_{ij} \otimes u^*_{ij} $
denote the conjugate representation of $u$.
 Then by \cite{free},
the interior tensor product representation
$$
u \otimes_{in} v \otimes_{in} \bar{u}
 = \sum_{ijklrs}  e^u_{ij} \otimes e^v_{kl} \otimes  e^{u}_{rs}
\otimes u_{ij} v_{kl} u^*_{rs}
$$
is irreducible.
%Let $h = h_1 * h_2$ be the Haar state of
%$G_1 \hat{*} G_2$, which is the free product of the Haar states
%$h_1, h_2$ of $G_1, G_2$.
Let $h_1 $ be the Haar state of $G_1 $.
Then
$$
h_1 \pi (u \otimes_{in} v \otimes_{in} \bar{u})
= \sum_{ijrs}  e^u_{ij} \otimes I_v \otimes  e^{u}_{rs}
 h_1(u_{ij} u^*_{rs}),
$$
where $I_v$ is the identity matrix
acting on the Hilbert space of $v$.
Since $u$ is of dimension greater than one, $u \otimes \bar{u}$
properly contains the trivial representation of $G_1$
(with multiplicity one). Hence
$$
h_1(u \otimes_{in} \bar{u}) =
\sum_{ijrs}  e^u_{ij} \otimes  e^{u}_{rs} h_1(u_{ij} u^*_{rs})
$$
is neither the identity matrix, nor the zero matrix
on $H_u \otimes H_{\bar{u}}$. (See also Theorem 5.7 of
Woronowicz \cite{Wor87b}.)
Therefore $h_1 \pi (u \otimes_{in} v \otimes_{in} \bar{u})$
is neither the identity matrix, nor the zero matrix
on $H_u \otimes H_v  \otimes H_{\bar{u}}$.
By \dfref{normal}, $G_1$ is not normal.

If $G_1$ has no non-trivial irreducible representations
of dimension greater than one, then by \cite{Wor87b}, $G_1$ is
the compact quantum group dual
of a discrete group $\Gamma$, i.e., $A_{G_1} = C^*(\Gamma)$.
By \cite{free}, every irreducible representation of $G$ is
of the form
$$w^{\lambda_1} \otimes w^{\lambda_2} \otimes \cdots
\otimes w^{\lambda_n},$$
where $w^{\lambda_i}$ is a  non-trivial representation belonging
to either the set $\Gamma $  or the set $ \hat{G}_2 $,
$w^{\lambda_i}$ and $w^{\lambda_{i+1}}$ being in different sets.
It is clear that
$$\pi(w^{\lambda_1} \otimes w^{\lambda_2} \otimes \cdots
\otimes w^{\lambda_n})$$
is a constant diagonal matrix, with the constant diagonal
entry equal to the product of those $w^{\lambda_i}$'s
that belong to $\Gamma$.
Use $[w^{\lambda_1} \otimes w^{\lambda_2} \otimes \cdots
\otimes w^{\lambda_n}]$ to denote this diagonal entry.
Since $h_{G_1}(\gamma) = 0$ if $\gamma$ is not
the neutral element of $\Gamma$,
%using the most general form of irreducible representations
%in \cite{W1} for $G= G_1 \hat{*} G_2$
one sees that \dfref{normal} is fulfilled.
That is,  $(G_1, \pi)$ is normal in $G= G_1 \hat{*} G_2$ and
${ A}_{G/{G_1}}$ is equal to the closure of the
linear span of entries of
the matrix $w^{\lambda_1} \otimes w^{\lambda_2} \otimes \cdots
\otimes w^{\lambda_n}$ such that
$[w^{\lambda_1} \otimes w^{\lambda_2} \otimes \cdots
\otimes w^{\lambda_n}]$ is the neutral element of $\Gamma$.
\qed

\medskip

The above result suggests that the free product
$G= G_1 \hat{*} G_2$ rarely has normal quantum groups.
This lead the author to conjecture that
free product of simple quantum groups is simple,
cf. Problem 4.3 of \cite{simpleprob}. 
Recently, Chirv\u{a}situ \cite{Chir12bpr} 
gave the following remarkable solution of this conjecture 
for simple compact quantum groups without center, generalizing his earlier results \cite{Chir12a} on simplicity of
the quotient quantum group of the universal unitary quantum group
by its center: 

\medskip
\noindent
{\em If $G_1$ and $G_2$ are simple compact quantum groups
with trivial center, then the free product $G_1 \hat{*} G_2$
is also simple.}

\medskip

In \cite{Chir12bpr}, Chirv\u{a}situ also proved that the 
quantum reflection group of Banica and Vergnioux \cite{B-Ver09} is simple. 
We refer the reader to his paper for details on these 
and other interesting results.

\medskip

Next we study normal quantum subgroups associated with crossed products.
Recall \cite{prod} that a {\bf discrete Woronowicz $C^*$-dynamical system }
is a triple
$(A, \Gamma, \alpha)$, where $A$ is a Woronowicz $C^*$-algebra,
 $\Gamma$ is a discrete group, and $\alpha$ is a homomorphism  from
$\Gamma$ to the automorphism group the Woronowicz $C^*$-algebra $A$,
i.e., automorphisms of the $C^*$-algebra $A$ that preserves the coproduct on $A$.
For such a dynamical system, it is shown in \cite{prod}
that the crossed product $C^*$-algebra $A \rtimes_\alpha \Gamma$ is also a
Woronowicz $C^*$-algebra, i.e., a compact quantum group, to abuse terminology.
The dense Hopf subalgebra of $A \rtimes_\alpha \Gamma$ is
${\mathcal A} \rtimes_\alpha \Gamma$.

\medskip

Similar to the second situation in
\propref{normalfreeproduct}, we have

\bgprop
\label{crossed-product}
Consider a crossed product $A \rtimes_\alpha \Gamma$
of a compact quantum group $A$ by a discrete group $\Gamma$.
Then $C^*(\Gamma)$ is a normal quantum subgroup
of $A \rtimes_\alpha \Gamma$
with quotient $A$ via the morphism
$$
\pi: A \rtimes_\alpha \Gamma \longrightarrow C^*(\Gamma),
\; \; \;
\pi(a \gamma) = \varepsilon(a) \gamma,
$$
where $\varepsilon$ is the counit of $A$ and $a \in A,
\; \gamma \in \Gamma$.
In the language of \cite{Andrus95a}, we have an
exact sequence of Hopf algebras
$$
{\mathbb C} \longrightarrow
{\mathcal A} \longrightarrow
{\mathcal A} \rtimes_\alpha \Gamma  \longrightarrow
{\mathbb C}  \Gamma \longrightarrow
{\mathbb C}.
$$

\ndprop

More generally, we have the following result that reduces
to the above proposition by taking $I=0$ and $K$ the one element group
\iffalse
,
bringing us more close to \propref{product-group}
\fi
(we use the notation in \thmref{Hopf $*$-ideals vs. Woronowicz $C^*$-ideals}
for correspondence of ideals).

\bgprop
Let $(A, \Gamma, \alpha)$ be a discrete Woronowicz $C^*$-dynamical system.
Let $I \lhd A$ be an $\alpha$-invariant Woronowicz $C^*$-ideal so that
$A/I$ is a normal quantum subgroup of $A$ with quotient quantum group $B$.
Let $K$ be a subgroup of the kernel of $\alpha$.
Let $\tilde \alpha$ be the evident action of
$\Gamma/K$ on $A/I$ obtained from $\alpha$.
Then $A/I \rtimes_{\tilde \alpha} \Gamma/K$
is a normal quantum subgroup of $A \rtimes_\alpha \Gamma$ with
quotient $B \rtimes_\alpha K$.
In the language of \cite{Andrus95a}, we have an exact sequence of Hopf algebras
$$
{\mathbb C} \longrightarrow
{\mathcal B} \rtimes_\alpha K \longrightarrow
{\mathcal A} \rtimes_\alpha \Gamma  \longrightarrow
{\mathcal A}/{\mathcal I} \rtimes_{\tilde \alpha} \Gamma/K \longrightarrow
{\mathbb C}.
$$
\ndprop
\noindent
\pf
Consider the morphisms
$$
\pi : A \longrightarrow A/I \rtimes_{\tilde \alpha} \Gamma/K ,
\; \; \text{and} \; \;
U : \Gamma \longrightarrow A/I \rtimes_{\tilde \alpha} \Gamma/K ,
$$
defined by $\pi(a) = {\tilde a}$, $U_\gamma = {\tilde \gamma}$.
Here
${\tilde a}$ is the image in $A/I$ of an element $a \in A$
viewed as an element of $A/I \rtimes_{\tilde \alpha} \Gamma/K$
via inclusion; %and
${\tilde \gamma}$ is the image in $\Gamma/K$ of
an element $\gamma \in \Gamma$ viewed as
an element of $A/I \rtimes_{\tilde \alpha} \Gamma/K$
via inclusion.
One can verify that $(\pi, U)$ is a covariant representation, i.e.,
$$U_\gamma \pi(a) U_\gamma^{-1} = \pi(\alpha_\gamma(a)).$$
Hence there is a surjection
$$
\pi \times U: A \rtimes_\alpha \Gamma
\longrightarrow A/I \rtimes_{\tilde \alpha} \Gamma/K ,
$$
extending $\pi$ and $U$.
This morphism preserves the coproducts because its restrictions
to $A$ and $C^*(\Gamma)$ do.
This shows that
$A/I \rtimes_{\tilde \alpha} \Gamma/K$
is a quantum subgroup of $A \rtimes_\alpha \Gamma$ (true under
only the assumption that $A/I $ is a quantum subgroup of $A$).

Let $A=A_G$ and assume that
$A/I = A_N$ is a normal quantum subgroup of $A$ with quotient $B = A_{G/N}$.
Let $\pi_0$ be the morphism from $A_G$ to $A_N$.
Let $\theta$ denote the surjection $\pi \times U$ found above.
Let $h = h_{A/I} \rtimes_{\tilde \alpha} h_{\Gamma/K}$ be the
Haar state on
$A/I \rtimes_{\tilde \alpha} \Gamma/K $ (cf. \cite{prod}).
Then every irreducible representation of the quantum group
$A \rtimes_{ \alpha} \Gamma $ is of the form
$(u^\lambda_{ij} {\gamma})$, where $(u^\lambda_{ij})$ is an irreducible
representation of dimension $d_\lambda$ of the quantum group $A$ and
${\gamma} \in \Gamma$ (cf. \cite{prod} as well).
Let $S(N)$ be the set of $\lambda$'s such that $h_N\pi_0(u^\lambda)$
is $I_{d_\lambda}$,
as in the proof of (4)$\Rightarrow$(3) in of Proposition 2.1 in \cite{simple}.
Then
$$\displaystyle{
(h \theta(u^\lambda_{ij} \gamma)) =
(h({\tilde u^\lambda_{ij}}) h({\tilde \gamma})) =
                    \begin{cases}
                 I_{d_\lambda} & \text{if \ $\lambda \in S(N) , \;
                                             \gamma \in K  ,$}\\
                 0             & \text{otherwise.}
                     \end{cases}}
$$
Hence by Proposition 2.1 in \cite{simple} and its proof which asserts
in part that
$$
 {\mathcal A}_{N \backslash G} = {\mathcal A}_{G/N}
= \bigoplus \{ {\mathbb C} u^\lambda_{ij} \;
| \; \lambda \in S(N), i, j = 1, \cdots, d_\lambda \} ,
$$
we conclude that
$A/I \rtimes_{\tilde \alpha} \Gamma/K$
is a normal quantum  subgroup of $A \rtimes_\alpha \Gamma$ with
quotient $B \rtimes_\alpha K$.

The last statement on exact sequence of Hopf algebras now follows from
\thmref{Hopf $*$-ideals vs. Woronowicz $C^*$-ideals} and remarks after
\lmref{reconstruct-N-from-G/N}.
\qed
\vv
Note that $A$ is not a normal quantum subgroup of $A \rtimes_\alpha \Gamma$
but rather a quotient quantum group, unlike the semi-direct product of groups.

\medspace
\medspace
\noindent
{\bf Acknowledgments.}
The author would like to record here his gratitude and indebtedness to
the anonymous referee who carefully read the manuscript and made
a number of very helpful suggestions and corrections that lead to
a substantial improvement of the paper.

The work reported here was started many years ago and most results were obtained
at the same time as those in \cite{simple}. It was
supported in part by the National Science Foundation grant DMS-0096136.
The author had opportunity to revisit this work
thanks to Hanfeng Li and Yi-Jun Yao who invited him 
to Chongqing University and Fudan University respectively during the summer of 2012.
The author would also like to thank Hanfeng Li and Dechao Zheng at
Chongqing University, and Yi-Jun Yao and Guoliang Yu at Fudan University for
their hospitality and for creating a congenial atmosphere.
He also wish to record his thanks to Huichi Huang, Pan Ma and
Qinggang Ren for making his visit enjoyable.
After an early version of this work was completed, Chirvasitu kindly informed
the author of his preprint \cite{Chir11pr}, in which he had proved the faithful
flatness conjecture stated in the early version of this paper as well
as in Conjecture 1 on p.3329 of \cite{simple}.

%\newpage


\begin{thebibliography}{99}



\bibitem{Andrus95a} Andruskiewitsch, N., Devoto, J.:
{Extensions of Hopf algebras},
{\em Algebra i Analiz}, {\bf 7}:1 (1995), 22--61.
translation in {\em St. Petersburg Math. J.} {\bf 7}:1 (1996),  17--52

\bibitem{B-Ver09}
Banica, T. and Vergnioux, R.:
Fusion rules for quantum reflection groups.
{\em J. Noncommut. Geom.} {\bf 3}  (2009), 327-359.

\bibitem{Bi03a} Bichon, J.:
Quantum automorphism groups of finite graphs.
{\em Proc. Amer. Math. Soc.} {\bf 131} (2003), no. 3, 665--673

\bibitem{Bi08a} Bichon, J.:
Algebraic quantum permutation groups.
{\em Asian-Eur. J. Math.} {\bf 1}  (2008),  no. 1, 1--13.

\bibitem{Bourbaki72} Bourbaki, N.:
{\em Commutative Algebra}, Springer, 1972


\bibitem{Chir11pr} Chirv\u{a}situ, A.:
Cosemisimple Hopf algebras are faithfully flat over Hopf subalgebras.
arXiv:1110.6701


\bibitem{Chir12a} Chirv\u{a}situ, A.:
    Free unitary groups are (almost) simple,
J. Math. Phys. 53, 123509 (2012)


\bibitem{Chir12bpr} Chirv\u{a}situ, A.:
Centers, cocenters and simple quantum groups, arXiv:1212.4763



\bibitem{Childs} Childs, Lindsay N.:
{\em Taming Wild Extensions: Hopf Algebras and
Local Galois Module Theory,}
American Mathematical Society,  2000.



\bibitem{CDPS12pr}
 Cirio, L.S., D'Andrea, A., Pinzari, C., Rossi, S.:
Connected components of compact matrix quantum groups and finiteness conditions.
arXiv:1210.1421


\bibitem{Dascalescu}
D\u{a}sc\u{a}lescu, S. and   N\u{a}st\u{a}sescu, C. and Raianu, \c{S}.:
{\em Hopf Algebras}, Marcel Dekker, Inc., New York, Basel 2001.


\bibitem{DemazureGabriel70} Demazure, M. and Gabriel, P.:
Groupes Alg\'ebriques I, North Holland, Amsterdam, 1970


\bibitem{FrSkTo12}
Franz, U.; Skalski, A.; Tomatsu, R.: 
On square roots of the Haar state on compact quantum groups. 
{\em J. Pure Appl. Algebra} 216 (2012), no. 10, 2079-2093. 


\bibitem{KlimykSchmudgen} Klimyk, A. U.  and Schm\"{u}dgen, K.:
{\em Quantum Groups and Their Representations},
Springer--Verlag, 1997.


\bibitem{Mackey76} Mackey, G. W.:
{\em The Theory of Unitary Group Representation,}
The University of Chicago Press, 1976


\bibitem{Montgomery93a} Montgomery, S.:
{\em Hopf Algebras and Their Actions on Rings,}
CBMS Regional Conf. Series, No 82,
American Mathematical Society, Providence, R.I., 1993



\iffalse
\bibitem{Nichols78a} Nichols, Warren D.
Quotients of Hopf algebras.
{\em Comm. Algebra} {\bf 6} (1978), no. 17, 1789--1800.
\fi

\bibitem{ParshallWang91a} Parshall, B. and Wang, J.:
{\rm Quantum linear groups,}
{\em Memoirs AMS} {\bf 439}, 1991.

\bibitem{Pod95a} Podles, P.:
{\rm Symmetries of quantum spaces. Subgroups and
quotient spaces of quantum $SU(2)$ and $SO(3)$ groups,}
{\em Commun. Math. Phys.} {\bf 170} (1995), 1-20.


\bibitem{Schauenburg00a} Schauenburg, Peter:
 Faithful flatness over Hopf subalgebras: counterexamples.
{\em Interactions between ring theory and representations of algebras}
(Murcia), 331--344, Lecture Notes in Pure and Appl. Math., 210,
Dekker, New York, 2000.

\bibitem{Schneid93a} Schneider, H.-J.:
Some remarks on exact sequences of quantum groups.
{\em Comm. Algebra} {\bf 21}:9 (1993), 3337--3357.

\bibitem{Sweedler} Sweedler, M. E.:
{\em Hopf Algebras,}
Benjamin, New York,  1969.

\bibitem{Takeuchi72} Takeuchi, Mitsuhiro:
A correspondence between Hopf ideals and sub-Hopf algebras,
{\em Manuscripta Math.} {\bf 7} (1972), 251-270.

\bibitem{Takeuchi79} Takeuchi, Mitsuhiro:
Relative Hopf modules-equivalences and freeness criteria,
{\em  J. Algebra} {\bf 60} (1979), no. 2, 452-471.


\bibitem{Takeuchi94a} Takeuchi, Mitsuhiro:
 Quotient spaces for Hopf algebras.
{\em Comm. Algebra} {\bf 22}:7 (1994), 2503--2523.

\bibitem{Daele95a} Van Daele, A.:
{\rm The Haar measure on a compact quantum group,}
{\em Proc. Amer. Math. Soc.} {\bf 123} (1995), 3125-3128.


\bibitem{uqg} Van Daele, A. and Wang, S. Z.:
{\rm Universal quantum groups,}
{\em International J. Math.} {\bf 7}:2 (1996), 255-264.

\bibitem{Wang} Wang, S. Z.:
{\em General Constructions of Compact Quantum Groups,}
Ph.D. Thesis, University of California at Berkeley, March, 1993.

\bibitem{free} Wang, S. Z.:
{\rm Free products of compact quantum groups,}
{\em Commun. Math. Phys.} {\bf 167}:3 (1995), 671-692.

\bibitem{prod} Wang, S. Z.:
{\rm Tensor products and crossed products of compact quantum groups,}
{\em Proc. London Math. Soc.} {\bf 71}:3 (1995), 695-720.


\bibitem{krein} Wang, S. Z.:
{\rm Krein duality for compact quantum groups,}
{\em J. Math. Phys.} {\bf 38}:1 (1997), 524-534


\bibitem{qsym} Wang, S. Z.:
{\rm Quantum symmetry groups of finite spaces},
math.OA/9807091,
{\em Commun. Math. Phys.} {\bf 195}:1 (1998), 195-211.

\bibitem{act} Wang, S. Z.:
{\rm Ergodic actions of universal quantum groups on operator algebras},
math.OA/9807093, {\it Commun. Math. Phys.} {\bf 203} (1999), 481-498.
%IHES, 1996. %Preprint, UC-Berkeley, 1997


\bibitem{simple} Wang, S. Z.:
{\rm  Simple Compact Quantum Groups I},
{\it J. Funct. Anal.} {\bf 256} (2009), no. 10, 3313--3341.

\bibitem{simpleprob} Wang, S. Z.:
 On the problem of classifying simple compact quantum groups,
{\it Banach Center Publ.} {\bf 98} (2012), 433-453


\bibitem{Wor87b} Woronowicz, S. L.:
{\rm Compact matrix pseudogroups,}
{\em Commun. Math. Phys.} {\bf 111} (1987), 613-665.


\bibitem{Wor98a} Woronowicz, S. L.:
{\rm Compact quantum groups,}
in {\em Quantum Symmetries},
Les Houches Summer School-1995, Session LXIV,
A. Connes, K. Gawedzki, J. Zinn-Justin, ed.,
Elsevier Science,
Amsterdam, 1998; pp845-884.





\end{thebibliography}
\end{document}